\documentclass[11pt]{article}
\usepackage{amssymb,amsmath,verbatim,amsthm}
\usepackage{tikz-cd}
\usepackage{tikz}
\usetikzlibrary{positioning}
\newtheorem{Theorem}{Theorem}[section]
\newtheorem{Lemma}[Theorem]{Lemma}
\newtheorem{Proposition}[Theorem]{Proposition}

\newtheorem{Corollary}[Theorem]{Corollary}

\newtheorem{Example}[Theorem]{Example}
\newtheorem{Definition}[Theorem]{Definition}

\newtheorem{Observation}[Theorem]{Observation}
\usepackage{enumerate}
\usepackage{color,xcolor}
\usepackage{indentfirst}
\usepackage{booktabs}
\usepackage{makecell}
\usepackage{array}
\usepackage{bm}
\newcommand{\ignore}[1]{}

\newcommand{\PG}{{\operatorname{PG}}}
\newcommand{\PGL}{{\operatorname{PGL}}}

\newcommand{\Aut}{{\operatorname{Aut}}}
\newcommand{\MAut}{{\operatorname{MAut}}}
\newcommand{\PAut}{{\operatorname{PAut}}}
\newcommand{\GDD}{{\operatorname{GDD}}}
\newcommand{\Stab}{{\operatorname{Stab}}}
\newcommand{\Tr}{{\operatorname{Tr}}}

\newcommand{\F}{\mathbb{F}}
\newcommand{\supp}{\operatorname{supp}}
\newcommand{\wt}{\operatorname{wt}}
\newcommand{\Aw}{\mathcal{A}}
\newcommand{\C}{\mathcal{C}}

\begin{document}
	\renewcommand{\theequation}{\thesection.\arabic{equation}}
	
	\title{\bf On $q$-Ary $t$-Designs from Codewords of Every Nonzero Weight in Linear Codes and Their Duals\footnote{Corresponding Author: Junling Zhou. Email: {\tt jlzhou@bjtu.edu.cn}}}
	\author{
		{ Xinghao Wu\footnote{School of Mathematics and Statistics, Beijing Jiaotong University,  Beijing  100044, China. (Email: {\tt 23111512@bjtu.edu.cn})}}, { Junling Zhou\footnote{School of Mathematics and Statistics, Beijing Jiaotong University,  Beijing  100044, China; Hebei Provincial Key Laboratory of Mathematical Theory and Analysis for Network and Data Science, Beijing Jiaotong University,  Beijing  100044, China
		}}\\
	}
	\date{ }
	\maketitle
	\begin{abstract}
		The paper investigates when the codewords of every nonzero weight in a linear code form a $q$-ary $t$-design and how to generalize the property to its dual code. We reduce the verification of this property to a prescribed range of low weights. Specifically, for MDS codes, we show that the $q$-ary $2$-design property at every nonzero weight is equivalent to that at any fixed nonzero weight. As an application, we prove that the codewords of every nonzero weight in a $[q+1,3,q-1]_q$ MDS code and its dual hold $q$-ary $2$-designs if and only if $q$ is even.	We also establish a duality principle. When $t\le\min\{d,d^{\perp}\}$, the codewords of every nonzero weight in a linear code form a $q$-ary $t$-design if and only if the same property holds for its dual code. We then investigate narrow-sense antiprimitive BCH codes $\C_{(q,q+1,\delta,1)}$ and their duals. For $q=2^m$ and $\delta=2^r$ with $1\le r<m$ and $\gcd(r,m)=1$, the codewords of every nonzero weight in these codes form $q$-ary $2$-designs. Furthermore, for $\delta=4$ and $q=2^{m}$ ($m\ge4$), we establish that the minimum-weight codewords of $\mathcal C_{(q,q+1,4,1)}$ form a $q$-ary $2$-design if and only if $m$ is odd.
		\medskip  
		
		\medskip
	\end{abstract}
	\section{Introduction}
	\setcounter{equation}{0}

	Combinatorial design theory and coding theory are closely related, with each providing useful constructions and methods for the other. On the one hand, incidence matrices of combinatorial designs can be used to construct linear codes over finite fields \cite{Key,Ding-DS,Tonchev1,Tonchev2}. On the other hand, the supports of codewords of a fixed weight in a code may form $t$-designs. The celebrated Assmus--Mattson theorem \cite{AM1969} provides sufficient conditions for such constructions. Since then, numerous families of linear and nonlinear codes supporting $t$-designs have been investigated \cite{Cameron,Ding-Steiner,Ding-five,Ding,Pless}. In recent years, significant progress has been made in constructing infinite families of $t$-designs from MDS, almost MDS, and near MDS codes \cite{DingTang2020,HengWang2023,XuCaoQu2022}. In particular, Tang and Ding \cite{BCH-4-design} constructed an infinite family of BCH codes supporting $4$-designs, while Yan and Zhou \cite{Yan} obtained further families of $3$-designs and $4$-designs. The role of highly transitive automorphism groups in constructing designs from linear codes has also been investigated \cite{Ding, 3-tran-BCH}.

	While ordinary $t$-designs arising from linear codes depend only on codeword supports, $q$-ary $t$-designs	also take the nonzero coordinate values into account.	This notion was introduced by Delsarte \cite{dsd's',Delsarte1973}.
	Every $q$-ary $t$-design induces an ordinary $t$-design	through its supports, but the converse need not hold.
	For $q=2$, the two notions coincide.
	Goethals and van Tilborg \cite{Tilborg1975,UP} established connections between $t$-regular codes and $q$-ary $t$-designs.
	Assmus, Goethals, and Mattson \cite{punctured-t} studied related generalized designs.
	Nevertheless, $q$-ary $t$-designs arising from nonbinary linear codes remain much less studied than ordinary designs.
	In our earlier work \cite{wu}, in order to obtain $q$-ary $t$-designs from linear codes, we developed the	Standard Criterion and the Puncturing--Shortening Criterion, together with an approach based on automorphism groups.
	These methods yielded families of $q$-ary $2$-designs from MDS codes, perfect codes, one- and two-weight codes, extremal self-dual codes, and related code	families.
	For codes not covered by the general criteria, the automorphism-group approach also gave $q$-ary $2$-designs in doubly-extended Reed--Solomon (DRS) codes and the duals of BCH codes $\C_{(q,q+1,4,1)}$. 
	Observe that these codes were shown to hold $q$-ary $2$-designs only at the two smallest nonzero weights and that the BCH codes themselves were not studied in terms of $q$-ary designs.
	This motivates the present study of $q$-ary designs at all nonzero weights in linear codes and their duals.

	In this paper, we reduce the verification of $q$-ary $t$-designs at all nonzero weights in a linear code to a prescribed range of low weights.
	For MDS codes, we prove that the codewords of every nonzero weight form $q$-ary $2$-designs if and only if the codewords of any fixed nonzero weight form a $q$-ary $2$-design.
	We also establish a duality principle when $t\le \min\{d,d^{\perp}\}$: the codewords of every nonzero weight in a linear code form a $q$-ary $t$-design if and only if the same property holds for its dual code.
	We apply these results to MDS codes and BCH codes. For	any $[q+1,3,q-1]_q$ MDS code and its dual, we prove that the codewords of every nonzero weight form $q$-ary $2$-designs if and	only if $q$ is even; in this case, the strength of the $q$-ary designs cannot be increased.
	We next study narrow-sense antiprimitive BCH codes of length $q+1$ and their dual codes. Let $q=2^m$	and $\delta=2^r$, where $1\le r<m$ and $\gcd(r,m)=1$. 
	We show that the codewords of every nonzero weight in both $\mathcal C_{(q,q+1,\delta,1)}$ and $\mathcal C_{\{1,\ldots,\delta-1\}}$ form $q$-ary $2$-designs.
	Furthermore, for $\delta=4$ and $m\ge4$, we establish a necessary and sufficient condition, which significantly improves the result of Theorems~5.12 and~5.13 in \cite{wu} .
	
	The remainder of this paper is organized as follows.
	Section~2 recalls basic concepts and results on linear codes and designs. 
	Section~3 gives criteria that reduce the verification of $q$-ary $t$-designs to smaller weights.
	Section~4 studies the design properties of a linear code and its dual. 
	Section~5 treats $[q+1,3,q-1]_q$ MDS codes and their duals.
	Section~6 investigates $q$-ary $2$-designs from BCH codes and their duals, including a complete characterization for the case \(\delta=4\).
	Section~7 summarizes the main results and presents open problems.

	\section{Preliminaries}
	\setcounter{equation}{0}
	\subsection{Linear codes}
	Let $q$ be a prime power and let $\mathbb{F}_{q}$ be the finite field of order $q$. For a positive integer $n$, set $[n] = \{1,2,\dots,n\}$ and denote by $\mathbb{F}_{q}^{n}$ the $n$-dimensional vector space over $\mathbb{F}_{q}$.	
	A vector $\bm{x}\in\mathbb{F}_{q}^{n}$ is written as $\bm{x} = (x_{1},x_{2},\dots,x_{n})$.
	The \emph{(Hamming) weight} $\mathrm{wt}(\bm{x})$ is the number of nonzero entries of $\bm{x}$, and its support is supp$(\bm{x})=\{i\in[n] : x_{i}\ne0\}$. The \emph{Hamming distance} between $\bm{x},\bm{y}\in\mathbb{F}_{q}^{n}$ is $ d(\bm{x},\bm{y})= |\{i\in [n]:x_{i}\ne y_{i}\}|$.
	A \emph{linear code} $\C$ of length $n$ over $\mathbb{F}_q$ is a subspace of $\mathbb{F}_{q}^{n}$.
	We write $[n,k]_q$ for a linear code over $\F_{q}$ of length $n$ and dimension $k$, and $[n,k,d]_{q}$ when its minimum distance is $d$.
	With respect to the standard inner product $\langle \bm{x}, \bm{y} \rangle = \sum_{i=1}^{n}x_{i}y_{i}$, the \emph{dual code} of $\C$ is
	\[
	\C^{\perp} = \{\bm{x}\in \mathbb{F}_{q}^{n}: \langle \bm{x}, \bm{y} \rangle=0 \text{ for all } \bm{y}\in \C \}.
	\]
	If $\C\subseteq\C^{\perp}$ the code is called \emph{self-orthogonal}, and if $\C=\C^{\perp}$ it is \emph{self-dual}.
	
	For a linear code $\C$ of length $n$, let $W(\cal C)$ denote the set of its distinct nonzero weights in $\C$. For each weight $w$, let $A_w$ denote the number of codewords of weight $w$. The \emph{weight enumerator} of $\cal C$ is $A(z)=\sum_{i=0}^n A_i z^i$.
	Let $s$ be the number of distinct nonzero weights of $\C$.
	The minimum distance of the dual code, written $d^\perp = d(\C^{\perp})$, is the \emph{dual distance} of $\C$, and the number $s^{\perp} = s(\C^{\perp})$ of nonzero weights in $\C^{\perp}$ is the \emph{external distance} of $\cal C$ \cite{dsd's'}.
	A linear code with parameters $[n,k,d]_{q}$ is called \emph{maximum distance separable} (MDS) if it meets the Singleton bound, i.e., $d=n-k+1$.

	Let $q$ be a prime power and let $n\ge 2$ be an integer such that $\gcd(n,q)=1$. Further let $m=\operatorname{ord}_n(q)$	and let $\alpha$ be a primitive $n$th root of unity in	$\mathbb{F}_{q^m}$. For each integer $i$, let $M_i(x)$ denote the minimal polynomial of $\alpha^i$ over	$\mathbb{F}_q$. Given integers $2\le\delta\le n$ and $0\le h\le n-1$, define
	\[
	g_{(q,n,\delta,h)}(x)=\operatorname{lcm}\left(M_h(x),M_{h+1}(x),\ldots,M_{h+\delta-2}(x)\right).
	\]
	The cyclic code of length $n$ over $\mathbb{F}_q$ with generator polynomial $g_{(q,n,\delta,h)}(x)$ is called a \emph{BCH code} with designed distance 	$\delta$ and is denoted by $\C_{(q,n,\delta,h)}$ \cite{fundamentals,cyclic-MDS}.
	When $h=1$, it is called a \emph{narrow-sense BCH code}. A BCH code is called \emph{primitive} if $n=q^m-1$ and \emph{antiprimitive} if $n=q^m+1$.
	By the BCH bound, the minimum distance $d$ of $\C_{(q,n,\delta,h)}$ satisfies $d\ge\delta$.

	\subsection{$t$-Designs and $q$-ary $t$-designs}
	
	A $t$-$(n,w,\lambda)$ \emph{design} is a pair $(\mathcal{P},\mathcal{B})$, where $\mathcal{P}$ is a set of $n$ \emph{points} and $\mathcal{B}$ is a family of $w$-subsets of $\mathcal{P}$ (called \emph{blocks}), such that every $t$-subset of $\mathcal{P}$ is contained in exactly $\lambda$ blocks. 
	For simplicity, we usually denote a $t$-design solely by its block set $\mathcal{B}$. 
	If all blocks in $\mathcal{B}$ are distinct, the design is called simple. 
	Denote by $\binom{\mathcal{P}}{w}$ the set of all $w$-subsets of $\cal P$.  A well-known property is that a $t$-$(n,w,\lambda)$ design is also an $i$-$(n,w,\lambda_i)$ design for every $0\le i\le t$, where
	\begin{align}\label{lambda-i}
		\lambda_i = \lambda \binom{n-i}{t-i}\bigg/\binom{w-i}{t-i}.
	\end{align}
	
	These concepts can be extended to vectors over finite fields, leading to the notion of $q$-ary designs (see \cite{punctured-t,dsd's',Delsarte1973,NewCR,wu}).
	To distinguish them from their classical counterparts, designs as defined above are called \emph{ordinary} or \emph{classical} $t$-designs.
	For vectors $\bm{a},\bm{b}\in\mathbb{F}_q^n$, we say that $\bm{a}$ \emph{covers} $\bm{b}$ if $a_i = b_i$ whenever $b_i \neq 0$.
	\begin{Definition}\label{def-qtdesign}
		Let $\mathcal{A}$ be a nonempty set of vectors (called blocks) in $\mathbb{F}_{q}^{n}$, each of  weight $w$.
		Then $\mathcal{A}$ is called a \emph{$q$-ary $t$-$(n,w,\lambda)$ design}, also denoted by $t$-$(n,w,\lambda)_{q}$ design, if for every vector $\bm{x} \in \mathbb{F}_{q}^{n}$ of weight $t$, there are exactly $\lambda$ blocks in $\mathcal{A}$ that cover $\bm{x}$.
		Equivalently, for every $\bm{x}\in\F_q^n$ with $\wt(\bm{x})=t$, there are exactly $\lambda$ vectors $\bm{a}\in\Aw$ such that $d(\bm{a},\bm{x})=w-t$.
		The parameters $t$ and $\lambda$ are referred to as strength and index respectively.
	\end{Definition}
	
	If $\mathcal{A}$ is a $q$-ary $t$-$(n,w,\lambda)$ design, then, for every integer $i$ with $0 \le i \le t$, $\mathcal{A}$ also forms a $q$-ary $i$-$(n, w, \lambda_i)$ design, where
	\begin{equation}\label{q-lambda-i}
		\lambda_{i} = \lambda(q-1)^{t-i} \binom{n-i}{t-i} \bigg/\binom{w-i}{t-i},\quad |\Aw|=\lambda_{0}=\lambda(q-1)^{t} \binom{n}{t}\bigg/\binom{w}{t}
	\end{equation}
	If $\mathcal{A}\subseteq \mathbb{F}_{q}^{n}$ is a $q$-ary $t$-$(n,w,\lambda)$ design, then the collection of its supports forms a classical $t$-$(n,w,(q-1)^{t}\lambda)$ design (possibly non-simple) \cite[Lemma 2.2]{punctured-t}. An immediate consequence is that for $q=2$ the two notions coincide; a binary $t$-design is equivalent to an ordinary $t$-design.
	It is readily verified that a $q$-ary $t$-($n,w,\lambda$) design is equivalent to a group divisible design $\GDD_{\lambda}(t,w,(q-1)n)$ of type $(q-1)^{n}$, see \cite[Proposition 2.7]{wu}.

	\subsection{Designs from linear codes}
	We are interested in linear codes whose codewords of a given weight form a design.
	For a linear code $\C$ with parameters $[n,k,d]_{q}$, let  $\mathcal{A}_w(\C)$ be the set of codewords of weight $w$; thus $|\mathcal{A}_w(\C)|=A_w$.
	If the set $\mathcal{A}_{w}(\C)$ itself is a $q$-ary $t$-$(n, w, \lambda)$ design, we say that the code \emph{holds} a $q$-ary $t$-design (at weight $w$).
	Index the codeword coordinates with an $n$-set $\mathcal{P}$. For any $w$ with $A_w \neq 0$, denote by $\mathcal{B}_w(\C)$ the set of supports of the codewords in $\mathcal{A}_{w}(\C)$. 
	When $(\mathcal{P}, \mathcal{B}_w(\C))$ is a $t$-$(n, w, \lambda)$ design, we say that $\C$ \emph{supports} a $t$-design (at weight $w$).

	A useful observation is that, within a certain range, every support corresponds to exactly $q-1$ codewords.
	\begin{Lemma}\cite[Lemma 4.25]{Ding} \label{repeat_q-1}
		Let $\C$ be an $[n,k,d]_{q}$ code and let $h$ be the largest integer not exceeding $n$ such that $h-\left \lfloor \frac{h+q-2}{q-1} \right \rfloor < d$. 
		Then for every weight $w$ with $d\le w\le h$, each support of a codeword of weight $w$ is the support of exactly $q-1$ codewords in $\mathcal{A}_w(\C)$.
	\end{Lemma}
	
	The classical Assmus--Mattson theorem \cite[Theorem 4.2]{AM1969} provides a fundamental connection between linear codes and ordinary designs.
	\begin{Theorem}\emph{(Assmus--Mattson)}\label{AM-Theorem}
		Let $\C$ be an $[n,k,d]_{q}$ linear code with dual distance $d^{\perp}$, and let $h$ be defined as in Lemma~\ref{repeat_q-1} using $d$. Define $h^{\perp}$ similarly using $d^{\perp}$.
		Write $\sum_{i=0}^{n}A_i z^i$ and $\sum_{i=0}^{n}A^{\perp}_i z^i$ for the weight enumerators of $\C$ and $\C^{\perp}$.
		For a positive integer $t<d$, let $m$ be the number of indices $i$ with $1\le i\le n-t$ and $A^{\perp}_i\neq0$. If $t\le d-m$, then
		\begin{itemize}
			\item [(1)] for every nonzero weight $w\in W(\C)$ satisfying $d\le w\le h$, $\mathcal{B}_{w}(\C)$ forms a $t$-design;
			\item [(2)] for every nonzero weight $w\in W(\C^{\perp})$ satisfying $d^{\perp}\le w\le \min\{n-t,h^{\perp}\}$,  $\mathcal{B}_{w}(\C^{\perp})$ forms a $t$-design.
		\end{itemize}
	\end{Theorem}
	
	\begin{Proposition}\cite[Proposition 2.10]{wu}
		Let $\C$ be an $[n,k,d]_{q}$ code and let $h$ be defined as in Lemma~\ref{repeat_q-1}.
		For any $w\in W(\C)$ with $d\le w\le h$, if $\mathcal{A}_{w}(\C)$ forms a $q$-ary $t$-$(n,w,\lambda)$ design, then $\mathcal{B}_{w}(\C)$ forms a  $t$-$(n,w,(q-1)^{t-1}\lambda)$ design.
	\end{Proposition}
	
	A $q$-ary analogue of the Assmus-Mattson theorem provides a direct condition under which the codewords of every nonzero weight form $q$-ary designs. We recall the following criterion from \cite[Theorem~2.6]{wu}.
	\begin{Theorem}\label{Standard-Criterion}\emph{(Standard Criterion)}
		Let $\C$ be a $q$-ary linear code with minimum distance $d$, $s$ nonzero weights,  dual distance $d^{\perp}$ and external distance $s^{\perp}$.
		Assume that $d>s^{\perp}$ or $d^{\perp}>s$, and set $t = \max\{d-s^{\perp}, d^{\perp}-s\}$.	Then 
		\begin{enumerate}
			\item [(1)] for every nonzero weight $w\in W(\C)$, $\mathcal{A}_{w} (\cal C)$ forms a $q$-ary $t$-design;
			\item [(2)] for every nonzero weight $w\in W(\C^{\perp})$,  $\mathcal{A}_{w} (\mathcal{C^{\perp}})$ forms a $q$-ary $t$-design.
		\end{enumerate}
	\end{Theorem}
	
	\subsection{Monomial equivalence and automorphism groups}
	Another way to obtain designs from codes exploits symmetries.
	We begin with the notion of monomial equivalence.
	A \emph{monomial} matrix over $\mathbb{F}_{q}$ is a square matrix having exactly one nonzero entry from $\mathbb{F}_{q}$ in each row and column. 
	Any monomial matrix $M$ can be written as $M=DP$ or $M=PD_{1}$, where $D,D_{1}$ are diagonal matrices and $P$ is a permutation matrix.
	Two linear codes $\C_{1},\C_{2}\subseteq\mathbb{F}_{q}^{n}$ are called \emph{monomially equivalent} if there exists a monomial matrix $M$ over $\mathbb{F}_{q}$ such that $\C_{2} = \C_{1}M$. Monomially equivalent codes share the same weight distribution.
	
	The monomial automorphism group of a linear code is built from monomial transformations that preserve the code.
	Let $\C$ be an $[n,k,d]_{q}$ code, and let $\mathcal{P}$ be the set of its coordinate positions.
	Write a codeword $\bm{c}\in\C$ as $\bm{c}=(c_{x})_{x\in\mathcal{P}}$.
	The \emph{permutation automorphism group} $\PAut(\C)$ is the subgroup of $\operatorname{Sym}(\mathcal{P})$ consisting of all permutations $g$ such that $(c_{g^{-1}x})_{x\in\mathcal{P}}\in\C$ for every $\bm{c}\in\C$.
	The \emph{monomial automorphism group} $\MAut(\C)$ is the subgroup of $(\mathbb{F}_{q}^{*})^{n}\rtimes\operatorname{Sym}(\mathcal{P})$ formed by all pairs $((a_{x})_{x\in\mathcal{P}};g)$ satisfying $(a_{x}c_{g^{-1}x})_{x\in\mathcal{P}}\in\C$ for all $\bm{c}\in\C$.
	Let $\operatorname{Gal}(\mathbb{F}_{q})$ denote the Galois group of $\mathbb{F}_{q}$ over its prime field.
	The \emph{automorphism group} $\operatorname{Aut}(\C)$ is the subgroup of $(\mathbb{F}_{q}^{*})^{n}\rtimes(\operatorname{Sym}(\mathcal{P})\times\operatorname{Gal}(\mathbb{F}_{q}))$ consisting of all elements $((a_{x})_{x\in\mathcal{P}};g,\gamma)$ such that $(a_{x}\,\gamma c_{g^{-1}x})_{x\in\mathcal{P}}\in\C$ for every $\bm{c}\in\C$.
	
	The group $\MAut(\C)$ (respectively, $\Aut(\C)$) is called \emph{$t$-transitive} if for any two ordered $t$-tuples  of distinct coordinates, there exists an element of $\MAut(\C)$ (respectively, $\Aut(\C)$) whose permutation part sends the first tuple  to the second.
	Clearly, if the automorphism group is $t$-transitive, then it is also $i$-transitive  for all $1\le i\le t$.
	A $1$-transitive group is simply called transitive.
	
	A classical result asserts that $t$-transitivity of the automorphism group yields ordinary $t$-designs.
	\begin{Theorem}\cite[Theorem 8.4.7]{fundamentals}\label{transitive-t-design}
		Let $\C\subseteq\mathbb{F}_{q}^{n}$ be a linear code with a $t$-transitive automorphism group. Then every nonempty set $\mathcal{B}_{w}(\C)$ with $w\ge t$ forms a $t$-design.
	\end{Theorem}
	For $q$-ary designs, $t$-transitivity of the code is not sufficient and a strengthened criterion was obtained in \cite{wu}.
	\begin{Theorem}\cite[Theorem 5.1]{wu}\label{transitive-design}
		Let $\C\subseteq\mathbb{F}_{q}^{n}$ be a linear code with a $t$-transitive automorphism group. 
		Fix a $t$-subset $T$ of coordinates. If for every vector $\bm{x}\in\mathbb{F}_{q}^{n}$ of weight $t$ with support $T$ there are exactly $\lambda$ codewords in $\mathcal{A}_{w}(\C)$ covering $\bm{x}$, then $\mathcal{A}_{w}(\C)$ forms a $q$-ary $t$-design.
	\end{Theorem}

	The following theorem establishes that monomial equivalence preserves the $q$-ary $t$-design property at each fixed weight. 
	
	\begin{Theorem}\label{monomial-qtD}
		Let $q$ be a prime power, and let $\C\subseteq\mathbb{F}_{q}^{n}$ be a linear code such that $\mathcal{A}_{w}(\C)$ forms a $q$-ary $t$-$(n,w,\lambda)$ design.
		If $\mathcal{D}\subseteq \mathbb{F}_{q}^{n}$ is a linear code monomially equivalent to $\C$, then $\mathcal{A}_{w}(\mathcal{D})$ also forms a $q$-ary $t$-$(n,w,\lambda)$ design.
	\end{Theorem}
	\begin{proof}	
		Since $\cal D$ is monomially equivalent to $\cal C$, there exists a monomial transformation $\phi$ mapping $\cal C$ onto $\cal D$, given by	$\phi(\bm{c}) = (\alpha_1 c_{\sigma(1)},\dots,\alpha_n c_{\sigma(n)})$, where $\sigma\in\operatorname{Sym}([n])$ and $\alpha_i\in\mathbb{F}_q^*$. 
		As a monomial transformation, $\phi$ preserves	both Hamming weights and the covering relation.
		Consequently, $\phi(\mathcal{A}_w(\C))=\mathcal{A}_w(\mathcal{D})$. 
		For any $\bm{y}\in\mathbb{F}_q^n$ with	$\mathrm{wt}(\bm{y})=t$, let $\bm{x}=\phi^{-1}(\bm{y})$.
		Since $\mathrm{wt}(\bm{x})=t$ and $\phi$ preserves the covering relation, we have
		\begin{align*}
			\bigl|\{\bm{d}\in\mathcal{A}_w(\mathcal{D}):\bm{d}\text{ covers }\bm{y}\}\bigr|
			=\bigl|\{\bm{c}\in\mathcal{A}_w(\C):\bm{c}\text{ covers }\bm{x}\}\bigr|=\lambda.
		\end{align*}
		Since $\bm{y}$ is arbitrary, $\mathcal{A}_w(\mathcal{D})$ forms a	$q$-ary $t$-$(n,w,\lambda)$ design.
	\end{proof}

	For convenience, we summarize some notations and definitions adopted throughout the paper in Table \ref{table-all}.
	\begin{table}[ht]
		\caption{Table of definitions and notations}\vspace{0.2cm}	
		\footnotesize
		\centering
		\renewcommand\arraystretch{1.2}
		\begin{tabular}{cc}
			\toprule
			Notation & Meaning   \\
			\midrule
			$\mathbb{F}_{q}$ & Finite field with $q$ elements   \\
			$\mathbb{F}_{q}^{*}$ & $\mathbb{F}_{q}\setminus\{0\}$  \\
			$\mathbb{F}_{q}^{n}$ & $n$-Dimensional vector space over $\mathbb{F}_{q}$ \\
			$\binom{\mathcal{P}}{k}$ & Set of all $k$-subsets of $\mathcal{P}$  \\
			$[n]$ & The set $\{1,2,\dots,n\}$  \\
			$\mathrm{wt}$($\bm{x}$) & Hamming weight of vector $\bm{x}$  \\
			$\mathrm{supp}(\bm{x})$ & Support of vector $\bm{x}$  \\
			$d=d(\C)$ &Minimum distance of $\C$  \\
			$d^{\perp}=d(\C^{\perp})$ & Dual distance of $\C$   \\
			$s=s(\C)$ & Number of nonzero weights in $\cal C$  \\
			$s^{\perp}=s(\C^{\perp})$ & External distance of $\C$ \\
			$A_{w}$ & Number of codewords of weight $w$ in $\cal C$  \\
			$\mathcal{A}_{w}(\C)$ & Set of all codewords  of weight $w$ in $\cal C$  \\
			$\mathcal{B}_{w}(\C)$ & Set of supports of all codewords in $\mathcal{A}_{w}(\C)$  \\
			$B_{\bm{x},i}$ & Number of codewords at distance $i$ from $\bm{x}$ \\
			$B_{\bm{x},i}(w)$ & Number of weight-$w$ codewords at distance $i$ from $\bm{x}$  \\
			$W(\C)$ & Set of nonzero weights in $\C$ \\
			$U$ & Subgroup of $\mathbb{F}_{q^{2}}^{*}$ comprised of all $(q+1)$-th roots of unity  \\	
			$\delta_{a,b}$ & Kronecker delta  \\
			$\psi$ & Additive character of $\F_{q}$\\
			\bottomrule
		\end{tabular}
		\label{table-all}
	\end{table}

	\section{Low-weight Reduction for $q$-Ary Designs at All Weights}
	\setcounter{equation}{0}
	The Standard Criterion provides an approach to obtaining $q$-ary $t$-designs at all weights of linear codes. The main objective of this section is to establish a reduction method of identifying $q$-ary designs at all weights, whenever the Standard Criterion does not take effect.
	We first  relate the $q$-ary design property to outer distributions and regularity of linear codes. Then, we use coset weight distributions and the external distance to reduce the verification of $q$-ary $t$-designs at all nonzero weights to a prescribed range of low weights. As a consequence, we derive a corresponding characterization for MDS codes.

	\subsection{$t$-Regular codes and $q$-ary $t$-designs}
	We now study the design property through the outer distribution of a code.
	Let $\C$ be a linear $[n,k,d]_{q}$ code.
	For $\bm{x}\in\mathbb{F}_q^n$, define the outer distribution by
	\begin{equation}\label{Bxi}
		B_{\bm{x},j} = |\{\bm{c} \in \C : d(\bm{c},\bm{x})=j \}|,\qquad 0\leq j\leq n.
	\end{equation}
	We further denote by $B_{\bm{x},j}(w)$ the number of codewords of weight $w$ at distance $j$ from $\bm{x}$.	
	Since the zero codeword contributes exactly when $j=\wt(\bm{x})$, we have
	\begin{equation}\label{Bxj}
		B_{\bm{x},j}= \delta_{j,\wt(\bm{x})} + \sum_{w\in W(\C)} B_{\bm{x},j}(w),
	\end{equation}
	where $\delta_{a,b}$ denotes the Kronecker delta.
	
	We recall the notion of $t$-regularity in terms of the outer distribution.
	\begin{Definition}
		Let $\C$ be a linear $[n,k,d]_q$ code with covering radius $\rho$, and let $t$ be an integer. 
		The code $\C$ is called \emph{$t$-regular} if, for each $0\le j\le n$, the value $B_{\bm{x},j}$ depends only on $d(\bm{x},\C)$ for all $\bm{x}\in\mathbb{F}_q^n$ with $d(\bm{x},\C)\le t$. 
		It is called \emph{completely regular} if it is $t$-regular for some $t\ge\rho$.
	\end{Definition}
	
	We state a local counting formula for a $q$-ary $t$-design. 
	
	\begin{Lemma}\cite[Lemma 3.10]{wu}	\label{lem:local-pattern}
		Let $\mathcal{A}\subseteq\F_q^n$ be a $q$-ary
		$t$-$(n,w,\lambda)$ design.  
		Let $A,B,C\subseteq[n]$ be pairwise disjoint subsets with $|A|=a,\ |B|=b,\ |C|=c$, where $ a+b+c\le t$.
		For each $i\in A\cup B$, fix $\alpha_i\in\F_q^*$. Then the number of vectors $\bm{v}\in\mathcal{A}$ satisfying 
		\[ 
		v_i=\alpha_i \quad (i\in A),\qquad v_i\in\F_q^\ast\setminus\{\alpha_i\}\quad (i\in B),\qquad v_i=0\quad (i\in C) 
		\]
		is independent of the particular choice of the $\alpha_i$s and depends only on $a,b,c$. Denote this number by $\lambda_{a,b,c}$. More precisely,
		\begin{align*}
			\lambda_{a,b,c}=\lambda(q-2)^b (q-1)^{t-a-b}\frac{\binom{n-a-b-c}{w-a-b}}{\binom{n-t}{w-t}}.
		\end{align*}
	\end{Lemma}
	
	The following lemma is distilled from the proof of \cite[Theorem 3.11]{wu} and is further refined here.
	\begin{Lemma}\label{lem:B_xj(w) independent x}
		Let $\C$ be an $[n,k,d]_q$ linear code. 
		Suppose that $w\in W(\C)$ and $\mathcal{A}_w(\C)$ forms a $q$-ary $t$-$(n,w,\lambda)$ design. 
		Let $\bm{x}\in\F_q^n$ with $\wt(\bm{x})\le t$.
		Then, if $j$ and $w$ are fixed, the quantity $B_{\bm{x},j}(w)$ depends only on $\wt(\bm{x})$.
	\end{Lemma}
	
	\begin{proof}
		Fix $\bm{x}\in\F_{q}^{n}$ and set $r =\wt(\bm{x})\le t$.  For each $\bm{c}\in\Aw_w(\cal C)$, partition $\supp(\bm{x})$ into
		\begin{align*}
			A&=\{i\in\supp(\bm{x}):c_i=x_i\},\\
			B&=\{i\in\supp(\bm{x}):c_i\ne 0,\ c_i\ne x_i\},\\
			C&=\{i\in\supp(\bm{x}):c_i=0\}.
		\end{align*}
		Write $a=|A|, b=|B|$, and $c=|C|$. Then $a+b+c=r$ and $d(\bm{x},\bm{c})=w-a+c$.
		For fixed $a,b,c$, there are $\frac{r!}{a!b!c!}$ possible partitions of $\supp(\bm{x})$ into $A,B,C$. By Lemma~\ref{lem:local-pattern}, each such partition is realized by exactly $\lambda_{a,b,c}$ codewords, independently of the particular vector $\bm{x}$. Consequently,
		\[
		B_{\bm{x},j}(w) = \sum_{\substack{a,b,c\ge 0 \\ a+b+c=r\\ w-a+c=j}}\lambda_{a,b,c}\frac{r!}{a!b!c!}.
		\]
		The right-hand side depends on $\bm{x}$ only through $r=\wt(\bm{x})$, as required.
	\end{proof}
	
	It was shown in \cite[Theorem~3.11]{wu} that, when $d\ge 2t$, the code $\mathcal C$ is $t$-regular if and only if $\mathcal A_w(\mathcal C)$ forms a $q$-ary $t$-design for every $w\in W(\mathcal C)$. 
	In this paper, we further consider the case $t\le d<2t$ and establish several general results that hold whenever $d\ge t$; these results will be used later.

	\begin{Lemma}\label{lem:all Aw q-ary t-design->t-regilar}
		Let $\cal C$ be an $[n,k,d]_q$ linear code and let  $t\leq d$ be a positive integer.
		Suppose that $\mathcal{A}_{w}(\C)$ forms a $q$-ary $t$-design for every nonzero weight $w\in W(\C)$.
		Then,  $\cal C$ is $t$-regular. 
	\end{Lemma}
	
	\begin{proof}
		Fix $\bm{x}\in\F_q^n$ with $\wt(\bm{x})\le t$ and fix any integer $j$ with $0\le j\le n$. By Lemma~\ref{lem:B_xj(w) independent x}, for every nonzero $w\in W(\C)$, the quantity $B_{\bm{x},j}(w)$ depends only on $\wt(\bm{x})$. It follows from \eqref{Bxj} that the same is true for $B_{\bm{x},j}$.
		
		Now let $\bm{y}\in\F_q^n$ satisfy $d(\bm{y},\C)=r\le t$, and choose $\bm{e}\in\C$ such that $d(\bm{y},\bm{e})=r$. Set $\bm{x}=\bm{y}-\bm{e}$ so that $\wt(\bm{x})=r$. Since $\C$ is linear, translation by $\bm{e}$ permutes its codewords. Hence
		\begin{align*}
			B_{\bm{y},j} &= |\{\bm{c}\in \C : d(\bm{y},\bm{c}) = j \}| \\
			&= |\{\bm{c}\in\C : d(\bm{y}-\bm{e},\bm{c}-\bm{e}) = j\}| \\
			&= |\{\bm{c}'\in\C : d(\bm{x},\bm{c}') = j\}| \\
			&= B_{\bm{x},j}.
		\end{align*}
		Since $B_{\bm{x},j}$ depends only on $\wt(\bm{x}) =r$, the value $B_{\bm{y},j}$ depends only on $r=d(\bm{y},\C)$. Therefore $\C$ is $t$-regular.
	\end{proof}
	
	The preceding lemma shows that, for $t\le d$, $t$-regularity is necessary for the $q$-ary $t$-design property at every nonzero weight. We next establish a sufficient condition for this property in terms of the outer distributions.
	
	\begin{Lemma}\label{lem: Btj->all Aw q-ary t-design}
		Let $\cal C$ be an $[n,k,d]_q$ linear code, and let $t\le d$ be a positive integer.
		Let $\bm{x}\in\F_{q}^{n}$ satisfy $\wt(\bm{x})=t$.
		Suppose that, for each fixed $j$, $B_{\bm{x},j}$ depends only on $t$.
		Then  $\mathcal{A}_{w}(\C)$ forms a $q$-ary $t$-design for every nonzero weight $w\in W(\C)$.
	\end{Lemma}
	\begin{proof}
		We proceed by induction on the nonzero weights of $\C$. Fix $\bm{x}\in\F_q^n$ with $\wt(\bm{x})=t$.
		
		First consider $w=d$. If $d(\bm{x},\bm{c})=d-t$ for some nonzero $\bm{c}\in\C$, then $\wt(\bm{c})\le d(\bm{c},\bm{x})+\wt(\bm{x})=d$.
		Hence
		\begin{align*}
			B_{\bm{x},d-t}(d)=B_{\bm{x},d-t}-\delta_{d,2t},
		\end{align*}
		which is independent of the particular choice of $\bm{x}$. Then $\mathcal{A}_d(\C)$ forms a $q$-ary $t$-design.
		
		Now let $w>d$ be a nonzero weight of $\C$, and assume that $\mathcal{A}_v(\C)$ is a $q$-ary $t$-design for every nonzero weight $v<w$. If $d(\bm{x},\bm{c})=w-t$, then  $\wt(\bm{c})\le w$, so no codeword of weight greater than $w$ contributes to	$B_{\bm{x},w-t}$. Thus
		$$	B_{\bm{x},w-t}(w) = B_{\bm{x},w-t}-\delta_{w,2t}-\sum_{v=d}^{w-1} B_{\bm{x},w-t}(v).$$
		By assumption, the first term on the right-hand side is independent of $\bm{x}$, while the induction hypothesis together with Lemma~\ref{lem:B_xj(w) independent x} shows that each term in the sum is also independent of $\bm{x}$. Hence $B_{\bm{x},w-t}(w)$ is independent of the particular choice of $\bm{x}$. In other words, every vector of weight $t$ is covered by the same number of codewords in $\mathcal{A}_w(\C)$. Therefore $\mathcal{A}_w(\C)$ is a $q$-ary $t$-design, completing the induction.
	\end{proof}
	
	It remains to determine when $t$-regularity implies the codewords of every nonzero weight form $q$-ary $t$-designs  in the range $t\le d<2t$. In this range, vectors of weight $t$ need not lie at the same distance from $\mathcal C$, so $t$-regularity does not directly give the uniformity required by Lemma~\ref{lem: Btj->all Aw q-ary t-design}. The following result shows that the minimum-weight design property provides the necessary additional condition.
	\begin{Theorem}\label{t-regular+d->all}
		Let $\C$ be an $[n,k,d]_{q}$ linear code and let $t$ be a positive integer satisfying $t\le d< 2t$. Suppose that $\C$ is $t$-regular and that $\Aw_{d}(\C)$ forms a $q$-ary $t$-design. Then every nonempty set $\Aw_{w}(\C)$ forms a $q$-ary $t$-design.
	\end{Theorem}
	\begin{proof}
		Let $\bm{x}\in\F_{q}^{n}$ with $\wt(\bm{x})=t$. 
		Since $\Aw_{d}(\C)$ forms a $q$-ary $t$-design, there exists a codeword $\bm{c}\in\Aw_{d}(\C)$ that covers $\bm{x}$. 
		Thus $d(\bm{x},\bm{c}) = d-t$ and hence $d(\bm{x},\C)\le d-t$.
		Obviously, $d(\bm{x},\bm{0})=t> d-t$.
		If $\bm{c}'\in\C\setminus\{\bm{0}\}$, then $\wt(\bm{c}')\ge d$, and therefore $d(\bm{x},\bm{c}')\ge \wt(\bm{c}')-\wt(\bm{x})\ge d-t$.
		Thus $d(\bm{x},\C)=d-t< t$.
		Since $\C$ is $t$-regular and all vectors of weight $t$ lie at the same distance $d-t$ from $\C$, it follows that, for every $0\le j\le n$, the value $B_{\bm{x},j}$ is independent of the particular choice of $\bm{x}$. Lemma~\ref{lem: Btj->all Aw q-ary t-design} now implies that every nonempty set $\Aw_{w}(\C)$ forms a $q$-ary $t$-design.
	\end{proof}

	The following example shows that the minimum-weight	design condition in Theorem~\ref{t-regular+d->all} cannot be omitted. 
	It gives a $2$-regular code $\C$	with minimum distance $d=3$ for which $\Aw_3(\mathcal C)$ does not form a $q$-ary $2$-design.
	
	\begin{Example}
		Let $\mathbb{F}_4=\{0,1,\alpha,\alpha^2\}$,	where $\alpha^2=\alpha+1$. 	Consider the linear code $\C$ over	$\mathbb{F}_4$ with parity-check matrix
		\[
		H=\begin{pmatrix} 1&0&1&1\\	0&1&1&\alpha \end{pmatrix}.
		\]
		Then $\C$ is a $[4,2,3]_4$ MDS code with weight enumerator $A(z)=1+12z^3+3z^4.$
		The coset weight distributions corresponding to	distances $0$, $1$, and $2$ from $\mathcal C$ are, respectively,
		\[
		(1,0,0,12,3),\	(0,1,3,7,5),\	(0,0,6,4,6).
		\]
		Hence $\C$ is $2$-regular. However, $\Aw_3(\C)$ does not form a $4$-ary $2$-design.
		Indeed, if $\Aw_3(\C)$ forms a $4$-ary $2$-$(4,3,\lambda)$ design, then \eqref{q-lambda-i} yields $\lambda=\frac{2}{3}$, a contradiction.
		Thus, $2$-regularity alone does not	guarantee that the minimum-weight codewords form a $q$-ary $2$-design when $d=3$.
	\end{Example}

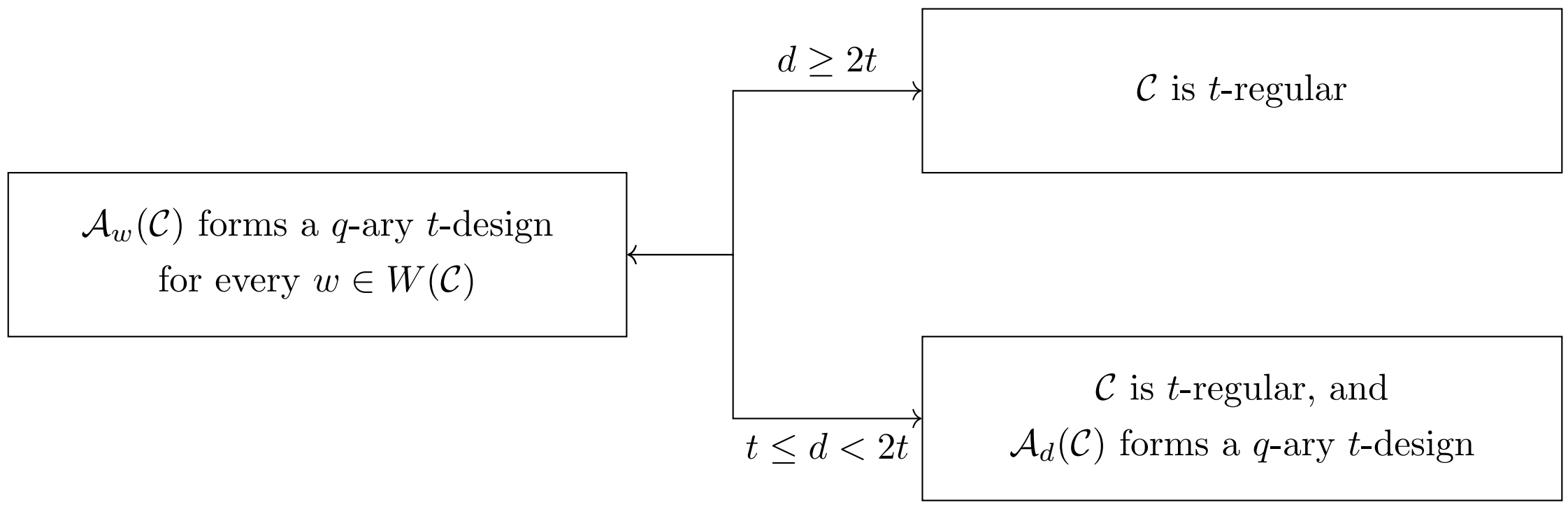
\begin{figure}[htbp]
		\centering
		\begin{tikzpicture}[
			font=\small,
			leftbox/.style={draw, align=center, text width=5.5cm,
				minimum height=1.55cm, inner sep=5pt},
			rightbox/.style={draw, align=center, text width=5.7cm,
				minimum height=1.55cm, inner sep=5pt}
			]
			\node[leftbox] (A) at (0,0) {\(\begin{gathered}
					\mathcal{A}_w(\mathcal C)\text{ forms a }q\text{-ary }t\text{-design}\\
					\text{for every }w\in W(\mathcal C)
				\end{gathered}\)};
			
			\node[rightbox] (B) at (8.75,1.55)
			{\(\mathcal C\text{ is }t\text{-regular}\)};
			
			\node[rightbox] (C) at (8.75,-1.55)
			{\(\begin{gathered}
					\mathcal C\text{ is }t\text{-regular, and}\\
					\mathcal A_d(\mathcal C)\text{ forms a }q\text{-ary }t\text{-design}
				\end{gathered}\)};
			
			\draw[<->] 
			(A.east) -- ++(1.0,0)
			|- node[pos=0.75, above] {\(d\ge 2t\)}
			(B.west);
			
			\draw[<->] 
			(A.east) -- ++(1.0,0)
			|- node[pos=0.75, below] {\(t\le d<2t\)}
			(C.west);
		\end{tikzpicture}
		\caption{Relationship between $t$-regularity and $q$-ary $t$-design property.}
		\label{fig:regularity-design-relations}
	\end{figure}

	Combining Lemma~\ref{lem:all Aw q-ary t-design->t-regilar}, Theorem~\ref{t-regular+d->all}, and \cite[Theorem~3.11]{wu}, we obtain the equivalences summarized in Fig.~\ref{fig:regularity-design-relations}.
	
	\subsection{Low-weight reduction}
	
	In this subsection, we  show that it suffices to consider a prescribed range of smaller weights to determine whether the codewords of every nonzero weight form  $q$-ary $t$-designs.
	
	\begin{Observation}\label{lem: Bxj=A(x+C)}
		Let $\cal C$ be an $[n,k,d]_q$ linear code. For any $\bm{x}\in\F_{q}^{n}$ and any integer $0\le j\le n$, the number of vectors of weight $j$ in the coset $\bm{x}+\C$ equals $B_{\bm{x},j}$.
	\end{Observation}
	\begin{proof}
		Since $\C$ is linear, $-\C=\C$. Hence
		\begin{align*}
			B_{\bm{x},j} &= |\{\bm{c} \in \C : d(\bm{x},\bm{c})=j \}|\\
			&= |\{\bm{c} \in \C : \wt(\bm{x}-\bm{c})=j \}| \\
			&= |\{\bm{v} \in \bm{x}+\C : \wt(\bm{v})=j \}| \\
			&=A_{j}(\bm{x}+\C).\qedhere
		\end{align*}
	\end{proof}

	The preceding observation allows us to use the following standard result for coset weight distributions.
	
	\begin{Lemma}\label{lem:coset weight s-1}\cite[Theorem~7.5.2]{fundamentals}
		Let $\mathcal C$ be an $[n,k]_q$ linear code with external distance $s^{\perp}$. For any $\bm{x}\in\mathbb{F}_q^n$, let $A_j(\bm{x}+\mathcal C)$ denote the number of vectors of weight $j$ in the coset $\bm{x}+\mathcal C$. Then the entire weight distribution of the coset is uniquely determined by $A_i(\bm{x}+\mathcal C)$ with $0\le i\le s^{\perp}-1$.
	\end{Lemma}

	\begin{Theorem}\label{th:s'+t-1}
		Let $\C$ be an $[n,k,d]_q$	linear code with external distance $s^{\perp}$. Let $t$ be a positive integer satisfying $d-s^{\perp}+1\le t\le d$.
		Suppose that, for every $w\in W(\C)$ satisfying $d\le w\le s^{\perp}+t-1$,	the set $\mathcal{A}_w(\C)$ forms a $q$-ary $t$-design.
		Then every nonempty set	$\mathcal{A}_w(\C)$ forms a $q$-ary $t$-design.
		Moreover, $\C$ is $t$-regular.
	\end{Theorem}
	
	\begin{proof}
		Fix an integer $r$ with $0\le r\le t$ and let $\bm{x}\in\mathbb{F}_q^n$	satisfy $\mathrm{wt}(\bm{x})=r$. If $\bm{c}\in\Aw_w(\cal C)$ satisfies $d(\bm{x},\bm{c}) = j$, then 
		\[
		w = \wt(\bm{c})\le d(\bm{x},\bm{c})+\wt(\bm{x}) =  j+r.
		\]
		Thus only codewords of weight at most $j+r$ can contribute to $B_{\bm{x},j}$, and therefore
		$$B_{\bm{x},j} = \delta_{r,j} + \sum_{w=d}^{j+r}B_{\bm{x},j}(w), \quad 0\le j\le n.$$
		Now let $0\le j\le s^{\perp}-1$. Every weight occurring in the sum above satisfies $d\le w\le j+r \le s^{\perp}+t-1$.
		By hypothesis, $\mathcal{A}_w(\C)$ forms a $q$-ary $t$-design for each such weight $w$. Since $r\le t$, Lemma~\ref{lem:B_xj(w) independent x}	implies that $B_{\bm{x},j}(w)$ depends on $\bm{x}$ only through $r$.
		Hence $B_{\bm{x},j}$ depends only on $r$ and $j$ for $0\le j\le s^{\perp}-1$.
		By Observation~\ref{lem: Bxj=A(x+C)}, we have $	B_{\bm{x},j}=A_j(\bm{x}+\C)$.	
		Thus, for every vector $\bm{x}$ of weight $r$,
		each $A_j(\bm{x}+\C)$, $0\le j\le s^{\perp}-1$,  is independent of the choice of $\bm{x}$.
		Lemma~\ref{lem:coset weight s-1} then implies that their entire weight distributions are identical.
		Therefore, for every $0\le j\le n$,	the quantity $B_{\bm{x},j}$ depends	only on $r=\mathrm{wt}(\bm{x})$
		whenever $r\le t$.
		
		In particular, for all vectors	$\bm{x}\in\mathbb{F}_q^n$ with	$\mathrm{wt}(\bm{x})=t$, the value $B_{\bm{x},j}$ is independent of the particular choice of $\bm{x}$ for every $0\le j\le n$.
		Since $t\le d$,	Lemma~\ref{lem: Btj->all Aw q-ary t-design}	shows that every nonempty set	$\mathcal{A}_w(\C)$ forms a $q$-ary $t$-design.
		Finally, Lemma~\ref{lem:all Aw q-ary t-design->t-regilar} yields that $\C$ is $t$-regular.
	\end{proof}

	In \cite[Proposition~6]{WuITW}, we proved that if the minimum-weight codewords of an MDS code form a $q$-ary $2$-design, then the codewords of every	nonzero weight also form $q$-ary $2$-designs.
	We now strengthen this result by showing that the $q$-ary $2$-design property at any fixed nonzero weight is equivalent to that at the minimum weight.

	\begin{Theorem}\label{pro:MDS-d->all}
		Let $\C$ be an $[n,k,d]_q$ MDS code with $d\ge 2$, and let $w\in W(\C)$ be an arbitrary  weight.
		Then $\Aw_w(\C)$ forms a	$q$-ary $2$-design if and only if $\Aw_d(\C)$ does.
	\end{Theorem}
	
	\begin{proof}
		If $k=1$, then $W(\C)=\{n\}$ and the assertion is immediate. Hence, assume that $k\ge 2$.
		
		Let $\bm{x}\in\F_q^n$ be an arbitrary vector of weight $2$, with $\supp(\bm{x})=\{i,j\}$, $x_i=a$, and $x_j=b$, where $a,b\in\F_q^*$.
		For each $w\in W(\C)$, define
		\[
		\lambda_{\bm{x}}(w)=|\{\bm{c}\in\Aw_w(\C): \bm{c}\text{ covers }\bm{x}\}|.
		\]
		For $\bm{c}\in\C$ and $0\le z\le k-1$, set
		\[
		N_{\bm{x}}(z)=|\{\bm{c}\in\C:\bm{c}\text{ covers }\bm{x},\ |\supp(\bm{c})|=n-z\}|.
		\]
		Clearly, $N_{\bm{x}}(z)=\lambda_{\bm{x}}(n-z)$.
		For each integer $0\le r\le k-1$, define
		\begin{align}\label{eq:Fx(r)1}
			F_{\bm{x}}(r)=\sum_{\substack{T\subseteq[n]\setminus\{i,j\}\\|T|=r}}
			\bigl|\{\bm{c}\in\C: c_i=a,c_j=b,c_\ell=0 \text{ for all }\ell\in T\}\bigr|.
		\end{align}
		Since every nonzero codeword of	$\C$ has at most $n-d=k-1$ zero coordinates,
		counting the pairs $(\bm{c},T)$ yields
		\begin{align}\label{eq:Fx(r)}
			F_{\bm{x}}(r)=\sum_{z=r}^{k-1}\binom{z}{r}N_{\bm{x}}(z).
		\end{align}		
		Applying the binomial inversion	formula to \eqref{eq:Fx(r)}, we obtain
		\[
		N_{\bm{x}}(z)=\sum_{r=z}^{k-1}(-1)^{r-z}\binom{r}{z}F_{\bm{x}}(r), \qquad 0\le z\le k-1.
		\]
		In particular, when $r=k-1$, \eqref{eq:Fx(r)} gives
		\begin{align*}
			F_{\bm{x}}(k-1)=N_{\bm{x}}(k-1)=\lambda_{\bm{x}}(d).
		\end{align*}
		
		For $0\le r\le k-2$, the prescribed	values in \eqref{eq:Fx(r)1} involve $r+2\le k$ distinct coordinates.
		Since $\C$ is MDS,	every assignment of values to any $r+2$ coordinates is realized	by exactly $q^{k-r-2}$ codewords.	Consequently,
		\[
		F_{\bm{x}}(r)=\binom{n-2}{r}q^{k-r-2},\qquad 0\le r\le k-2.
		\]
		Now fix $w\in W(\C)$ and put $\ell=n-w$.
		Since $d\le w\le n$, we have $0\le\ell\le k-1$. It follows that
		\begin{align*}
			\lambda_{\bm{x}}(w)=N_{\bm{x}}(\ell)=\sum_{r=\ell}^{k-2}(-1)^{r-\ell}\binom{r}{\ell}\binom{n-2}{r}q^{k-r-2}+(-1)^{k-1-\ell}\binom{k-1}{\ell}\lambda_{\bm{x}}(d),
		\end{align*}
		where the sum is understood	to be zero if $\ell=k-1$.
		The first term in this expression is independent of $\bm{x}$.
		Hence, for any two vectors $\bm{x},\bm{y}\in\F_q^n$	of weight $2$,
		\begin{align*}
			\lambda_{\bm{x}}(w)-\lambda_{\bm{y}}(w)=(-1)^{w-d}\binom{k-1}{n-w}\bigl(\lambda_{\bm{x}}(d)-\lambda_{\bm{y}}(d)\bigr),
		\end{align*}
		where we have used	$d=n-k+1$.
		
		Since $0\le n-w\le k-1$, the binomial coefficient $\binom{k-1}{n-w}$ is nonzero.
		Therefore, $\lambda_{\bm{x}}(w)$ is	independent of the choice of $\bm{x}$ with $\wt(\bm{x})=2$
		if and only if $\lambda_{\bm{x}}(d)$ is	independent of that choice.
		By the definition of a $q$-ary $2$-design, $\Aw_w(\C)$ forms a $q$-ary $2$-design if and
		only if $\Aw_d(\C)$ does.
	\end{proof}

	\begin{Corollary}\label{cor:MDS-some-all}
		Let $\C$ be an $[n,k,d]_q$ MDS code with $d\ge 2$. Then the following statements are equivalent:
		\begin{enumerate}[(1)]
			\item $\Aw_w(\C)$ forms a $q$-ary	$2$-design for some $w\in W(\C)$;
			\item $\Aw_w(\C)$ forms a $q$-ary	$2$-design for every $w\in W(\C)$.
		\end{enumerate}
	\end{Corollary}
	\begin{proof}
		Suppose that $\Aw_{w}(\C)$ forms a $q$-ary $2$-design for some $w\in W(\C)$.
		By Theorem~\ref{pro:MDS-d->all}, $\Aw_d(\C)$ also forms a $q$-ary $2$-design.
		Applying the same theorem again shows that	$\Aw_w(\C)$ forms a $q$-ary $2$-design for every $w\in W(\C)$.
		The converse is immediate.
	\end{proof}

	\section{$q$-Ary $t$-Designs in $\C$ and $\C^{\perp}$}
	\setcounter{equation}{0}
	The preceding section provides sufficient conditions under which the codewords of every nonzero weight form $q$-ary $t$-designs. A natural question is whether this property is preserved under duality. In the binary case,  if the codewords of every nonzero weight support a  $t$-design when $t\le d$, then the same holds for the dual code \cite[Theorem~13, p.~165]{cyclic-MDS}. In this section, we investigate $q$-ary $t$-designs arising	from a linear code and its dual. We use additive characters and Krawtchouk polynomials to show that, when $t\le\min\{d,d^{\perp}\}$, the codewords of every nonzero weight form $q$-ary $t$-designs in a linear code if and only if the same holds for its dual.

	Let $q=p^{e}$, where $p$ is prime. An \emph{additive character} of $\F_q$ is a group homomorphism $$\psi:\F_q\longrightarrow \mathbb{C}^{*},$$
	satisfying
	\begin{align}\label{psi(x+y)=p(x)p(y)}
		\psi(x+y)=\psi(x)\psi(y),\qquad x,y\in\F_q.
	\end{align}
	Let $\Tr_{q/p}$ denote the absolute trace from $\F_q$ to $\F_p$ and denote $\zeta_p=e^{2\pi\sqrt{-1}/p}$.
	The canonical additive character of $\F_q$ is $\psi_1(x)=\zeta_p^{\Tr_{q/p}(x)}$.
	Every additive character of $\F_q$ has the form	$$\psi_b(x)=\psi_1(bx)=\zeta_p^{\Tr_{q/p}(bx)}$$ for a unique $b\in\F_q$ \cite[Theorem~5.7]{Finite-Fields}. 
	The character $\psi_0$ is trivial, whereas $\psi_b$ is nontrivial for every $b\in\F_q^{*}$. In particular, $\psi(0)=1$ for every additive character $\psi$.
	
	The following standard orthogonality relation will be used throughout this section.
	
	\begin{Lemma}\cite[Theorem~5.4]{Finite-Fields}\label{add-char-sum}
		Let $\psi$ be a nontrivial additive character of $\F_q$. Then, for every $\alpha\in\F_q$,
		\begin{align*}
			\sum_{x\in\F_q}\psi(\alpha x)=
			\begin{cases}
				q, & \alpha=0,\\
				0, & \alpha\neq0.
			\end{cases}
		\end{align*}
		Consequently, if $\alpha\neq0$, then $\sum_{x\in\F_q^{*}}\psi(\alpha x)=-1$.
	\end{Lemma}

	We next calculate the fixed-weight character sum that occurs in the subsequent design-counting argument. Recall the $q$-ary Krawtchouk polynomial 
	$$
	K_r^{(m,q)}(s)=\sum_{\ell=0}^{r}(-1)^{\ell}(q-1)^{r-\ell}\binom{s}{\ell}\binom{m-s}{r-\ell}, \quad 0 \le s \le m.
	$$
	
	The following lemma is a direct generalization of \cite[Lemma 2.6.2]{fundamentals}.
	\begin{Lemma}\label{lem:support-character-sum}
		Let $\psi$ be a nontrivial additive character of $\F_q$, let $E\subseteq[n]$ with $|E|=m$, and let $\bm{z}\in\F_q^n$. Set $s=|E\cap\supp(\bm{z})|$.
		Then, for every $0\le r\le m$,
		\[
		\sum_{\substack{\bm{x}\in\F_q^n\\ \supp(\bm{x})\subseteq E\\ \wt(\bm{x})=r}}
		\psi\bigl(\langle\bm{z},\bm{x}\rangle\bigr)	=K_r^{(m,q)}(s).
		\]
	\end{Lemma}
	\begin{proof}
		Partition the vectors $\bm{x}$ according to their support $R\subseteq E$, $|R|=r$, and put
		\[
		\ell=|R\cap\supp(\bm{z})|.
		\]
		For a fixed support $R$, the nonzero coordinates of $\bm{x}$ vary independently over $\F_q^{*}$, and hence
		\[
		\sum_{\substack{\bm{x}\in\F_q^n\\ \supp(\bm{x})=R}}
		\psi\bigl(\langle\bm{z},\bm{x}\rangle\bigr)=\prod_{i\in R}\left(\sum_{x_i\in\F_q^{*}}\psi(z_i x_i)\right).
		\]
		By Lemma~\ref{add-char-sum}, a factor equals $-1$ when $z_i\neq0$ and equals $q-1$ when $z_i=0$. Therefore every support $R$ with $|R\cap\supp(\bm{z})|=\ell$ contributes $(-1)^{\ell}(q-1)^{r-\ell}$.
		There are $\binom{s}{\ell}\binom{m-s}{r-\ell}$ such supports. Summing over $\ell$ gives
		\[
		\sum_{\ell=0}^{r}(-1)^{\ell}(q-1)^{r-\ell}\binom{s}{\ell}\binom{m-s}{r-\ell}=K_r^{(m,q)}(s),
		\]
		as required.
	\end{proof}

	We now apply Lemma~\ref{lem:support-character-sum} to the codewords of fixed weights of the dual code. Throughout the remainder of this section, let $\psi$ be a fixed nontrivial additive character of $\F_q$.
	
	\begin{Lemma}\label{lem:dual-shell-fourier}
		Let $\C$ be an $[n,k,d]_q$ linear code, and let $t\le d$ be a positive integer. Suppose that $\Aw_w(\C)$ forms a $q$-ary $t$-design for every nonzero weight $w\in W(\C)$. For $0\le u\le n$ and $\bm{x}\in\F_q^n$, define
		\[
		\Phi_u(\bm{x})=\sum_{\bm{v}\in\Aw_u(\C^{\perp})}\psi\bigl(\langle\bm{x},\bm{v}\rangle\bigr).
		\]
		If $\wt(\bm{x})\le t$, then $\Phi_u(\bm{x})$ depends only on $\wt(\bm{x})$ and $u$.
	\end{Lemma}
	\begin{proof}
		For $\bm{v}\in\F_q^n$, additive-character orthogonality gives
		\begin{align}
			\frac{1}{|\C|}	\sum_{\bm{c}\in\C}	\psi(-\langle\bm{c},\bm{v}\rangle)=
			\begin{cases}
				1, & \bm{v}\in\C^{\perp},\\
				0, & \bm{v}\notin\C^{\perp}.
			\end{cases}
		\end{align}
		Indeed, if $\bm{v}\in\C^{\perp}$, then $-\langle\bm{c},\bm{v}\rangle=0$ for all $\bm{c}\in\C$, so every term $	\psi(-\langle\bm{c},\bm{v}\rangle) $ equals $1$. If $\bm{v}\notin\C^{\perp}$, then the map
		$\bm{c}\longmapsto-\langle\bm{c},\bm{v}\rangle$
		is a nonzero linear functional from $\C$ onto $\F_q$. Hence every element of $\F_q$ is attained exactly $|\C|/q$ times.  Therefore, by Lemma~\ref{add-char-sum},
		\[
		\sum_{\bm{c}\in\C}\psi(-\langle\bm{c},\bm{v}\rangle)=\frac{|\C|}{q}\sum_{a\in\F_{q}}\psi(a) = 0.
		\]
		Then applying Lemma~\ref{lem:support-character-sum} with $E=[n]$, we obtain
		\begin{align*}
			\Phi_u(\bm{x})&=
			\sum_{\substack{\bm{v}\in\F_q^n\\\wt(\bm{v})=u}}\left(\frac{1}{|\C|}\sum_{\bm{c}\in\C}\psi(-\langle\bm{c},\bm{v}\rangle)\right)
			\psi\bigl(\langle\bm{x},\bm{v}\rangle\bigr)	\notag\\
			&=\frac{1}{|\C|}\sum_{\bm{c}\in\C}\sum_{\substack{\bm{v}\in\F_q^n\\\wt(\bm{v})=u}}\psi(\langle\bm{x}-\bm{c},\bm{v}\rangle)\notag\\
			&=\frac{1}{|\C|}\sum_{\bm{c}\in\C}K_u^{(n,q)}\bigl(\wt(\bm{x}-\bm{c})\bigr)\notag\\
			&=\frac{1}{|\C|}\sum_{j=0}^{n}B_{\bm{x},j}K_u^{(n,q)}(j).
		\end{align*}
		
		By Lemma~\ref{lem:B_xj(w) independent x}, the hypothesis implies that whenever $\wt(\bm{x})\le t$, the numbers $B_{\bm{x},j}$ depend only on $\operatorname{wt}(\bm{x})$ and $j$.
		Therefore, the value $\Phi_u(\bm{x})$ depends only on $\wt(\bm{x})$ and $u$.
	\end{proof}
	
	Equipped with these lemmas, we now present a main result.
	
	\begin{Proposition}\label{th:dual-shell-design}
		Let $\C$ be an $[n,k,d]_{q}$ linear code, and let $t\le d$ be a positive integer. Suppose that $\Aw_{w}(\C)$ forms a $q$-ary $t$-design for every nonzero weight $w\in W(\cal C)$. Then, for every $u\in W(\mathcal{C^{\perp}})$ with $u\ge t$, $\Aw_{u}(\mathcal{C^{\perp}})$ forms a $q$-ary $t$-design.
	\end{Proposition}
	\begin{proof}
		Fix $u\in W(\C^\perp)$ with $ u\ge t$ and let $\bm{y}\in\F_q^n$ be any vector of weight $t$. Set $T=\supp(\bm{y})$. Define
		\[
		\lambda _{u,\bm{y}}	=	|\{\bm{v}\in\Aw_u(\C^\perp):	\bm{v}\text{ covers }\bm{y}\}|,
		\]
		For $\bm{v}\in\F_{q}^{n}$, coordinatewise additive-character orthogonality yields
		\begin{align*}
			\delta_{(\bm{v}|_T,\bm{y}|_T)}&=
			q^{-t}\prod_{i\in T}\left(\sum_{x_i\in\F_q}	\psi\bigl(x_i(v_i-y_i)\bigr)\right)\notag\\
			&=q^{-t}\sum_{\substack{\bm{x}\in\F_q^n\\ \supp(\bm{x})\subseteq T}}\psi\left(\sum_{i\in T}x_i(v_i-y_i)\right),
		\end{align*}
		where $\bm{v}|_{T}$ denotes the restriction of $\bm{v}$ to the coordinates in $T$.
		Substituting this identity into $\lambda _{u,\bm{y}}$ and interchanging the finite sums gives
		\begin{align}
			\lambda_{u,\bm{y}}&=q^{-t}\sum_{\substack{\bm{x}\in\F_q^n\\ \supp(\bm{x})\subseteq T}}\sum_{\bm{v}\in\Aw_u(\C^{\perp})}\psi\left(\sum_{i\in T}x_i v_i -  \sum_{i\in T}x_i y_i\right)\notag\\
			&=q^{-t}\sum_{\substack{\bm{x}\in\F_q^n\\ \supp(\bm{x})\subseteq T}}\psi\left(-\langle\bm{x},\bm{y}\rangle\right)\sum_{\bm{v}\in\Aw_u(\C^{\perp})}\psi(\langle\bm{x},\bm{v}\rangle)\notag\\
			&=q^{-t}\sum_{\substack{\bm{x}\in\F_q^n\\ \supp(\bm{x})\subseteq T}}\psi\left(-\langle\bm{x},\bm{y}\rangle\right)\Phi_u(\bm{x}).
			\label{eq:lambda-Phi}
		\end{align}
		Every vector $\bm{x}$ occurring in \eqref{eq:lambda-Phi} has weight at most $t$. By Lemma~\ref{lem:dual-shell-fourier}, $\Phi_u(\bm{x})$ depends only on $r=\wt(\bm{x})$ and $u$; denote this common value by $\phi_{u,r}$. Grouping the terms in \eqref{eq:lambda-Phi} according to $r$ yields
		\begin{equation}\label{eq:lambda-group-r}
			\lambda_{u,\bm{y}}=q^{-t}\sum_{r=0}^{t}\phi_{u,r}\sum_{\substack{\bm{x}\in\F_q^n\\ \supp(\bm{x})\subseteq T\\ \wt(\bm{x})=r}}
			\psi\left(-\langle\bm{x},\bm{y}\rangle\right).
		\end{equation}
		Since $T=\supp(\bm{y})$, Lemma~\ref{lem:support-character-sum}, applied with $E=T$ and $\bm{z}=-\bm{y}$, gives
		\[
		\sum_{\substack{\bm{x}\in\F_q^n\\\supp(\bm{x})\subseteq T\\\wt(\bm{x})=r}}\psi\bigl(-\langle\bm{x},\bm{y}\rangle\bigr)
		=K_r^{(t,q)}(t)=(-1)^r\binom{t}{r}.
		\]
		Therefore
		\[
		\lambda_{u,\bm{y}}=q^{-t}\sum_{r=0}^{t}(-1)^r\binom{t}{r}\phi_{u,r},
		\]
		which is independent of the particular vector $\bm{y}$ of weight $t$. Thus every vector of weight $t$ is covered by the same number of vectors in $\Aw_u(\C^{\perp})$, and hence $\Aw_u(\C^{\perp})$ forms a $q$-ary $t$-design.
	\end{proof}
	
	\begin{Theorem}\label{CCperp-qt}
		Let $\C$ be an $[n,k,d]_{q}$ linear code, and let $d^{\perp}$ be the minimum distance of $\C^{\perp}$. Let $t$ be a positive integer with $t\leq\min\{d,d^{\perp}\}$.
		Then every nonempty set $\Aw_{w}(\C)$ forms a $q$-ary $t$-design if and only if every nonempty set $\Aw_{u}(\C^{\perp})$ forms a $q$-ary $t$-design.
	\end{Theorem}
	\begin{proof}
		Suppose first that the codewords of every nonzero weight in $\C$ form a $q$-ary $t$-design. Since $d^{\perp}\ge t$, every nonzero weight $u$ of $\C^{\perp}$ satisfies $u\ge d^{\perp}\ge t$, and Proposition~\ref{th:dual-shell-design} implies that the codewords of every nonzero weight in $\C^{\perp}$ form a $q$-ary $t$-design.
		Conversely, apply the same argument to $\C^{\perp}$. Since $(\C^{\perp})^{\perp}=\C$ and $d\ge t$, Proposition~\ref{th:dual-shell-design} yields the converse implication.
	\end{proof}

	\section{$q$-Ary $2$-Designs from $[q+1,3]_{q}$ MDS Codes}
	\setcounter{equation}{0}
	
	In \cite[Theorem~4.13]{wu}, the relationship between hyperovals in $\PG(2,q)$ and $q$-ary $2$-designs was established for the associated $[q+2,3,q]_q$ MDS codes. In this section, we investigate the corresponding problem for $[q+1,3,q-1]_q$ MDS codes, whose generator matrices give rise to ovals in $\PG(2,q)$.
	We show that the codewords of every nonzero weight form $q$-ary $2$-designs if and only	if $q$ is even. We then extend this characterization to the	dual MDS codes. We also show that the codewords of any nonzero weight in these $[q+1,3,q-1]_q$ codes do not form a $q$-ary $3$-design.
	
	Recall that the points of a projective plane $\PG(2,q)$ are the one-dimensional subspaces of $\mathbb{F}_{q}^{3}$, and the lines are its two-dimensional subspaces. An \emph{oval} $\mathcal{O}$ is a set of $q+1$ points in $\PG(2,q)$, no three of which are collinear, whereas a \emph{hyperoval} $\mathcal{H}$ is a set of $q+2$ points with the same property \cite{two,Ding,PG-BOOK}. If $q$ is even,  every oval $\mathcal{O}$ has a unique \emph{nucleus} $K$: the $q+1$ tangent lines to $\mathcal{O}$ are concurrent at $K$, and $\mathcal{H}=\mathcal{O}\cup\{K\}$
	is a hyperoval. Conversely, every line of $\PG(2,q)$ meets a hyperoval in either $0$ or $2$ points. In particular, a line through the nucleus $K$ meets $\mathcal{O}$ in exactly one point and is therefore tangent to $\mathcal{O}$.
	
	Let $\mathcal{O}=\{P_1,\ldots,P_{q+1}\}$ be an oval. For each $P_i$, choose a nonzero column vector $\bm{p}_i\in\F_q^3$ representing $P_i=[\bm{p}_i]$. Define the associated code $\C_{\mathcal{O}}$ by the generator matrix
	\[
	G_{\mathcal{O}}=[\bm{p}_1\ \bm{p}_2\ \cdots\ \bm{p}_{q+1}].
	\]
	The following standard observation connects codeword weights with intersections in the projective plane.
	
	\begin{Lemma}\label{Ha}\cite[Theorem~2.43]{Ding}
		Let $\mathcal{S}=\{P_{1},\ldots,P_{n}\}$ be a set of points in $\PG(2,q)$ and let $G=[\bm{p}_1 \cdots \bm{p}_{n}]$ consist of fixed column representatives of these points. For $\bm{x}\in\mathbb{F}_{q}^{3}\setminus\{\bm{0}\}$, let
		\[
		\ell_{\bm{x}}=\{[\bm{y}]\in\PG(2,q):\langle \bm{x}, \bm{y} \rangle=0\}.
		\]
		If $\bm{c}=\bm{x}G$, then $\wt(\bm{c})=n-|\ell_{\bm{x}}\cap\mathcal{S}|$.
	\end{Lemma}
	Since no three points of an oval are collinear, any three columns of  $G_{\mathcal{O}}$ are linearly independent.
	Hence $\C_{\mathcal{O}}$ has parameters $[q+1,3,q-1]_q$. Conversely, the columns of a generator matrix of any $[q+1,3,q-1]_q$ MDS code give rise to an oval. Thus, every such code can be represented by an oval, with different choices of column representatives corresponding to monomially equivalent codes; see also \cite[Corollary~2.52]{PG-BOOK}. The weight distribution of such a code is determined by the MDS parameters \cite[Theorem 7.4.1]{fundamentals}.
	
	\begin{Lemma}\label{lem:oval-weight-distribution}
		Let $q\ge 3$ be a prime power and let $\C$ be a $[q+1,3,q-1]_{q}$ MDS code. Then the weight enumerator of $\C$ is 
		\begin{align}\label{eq:oval-weight-distribution}
			A(z)=1+\frac{q(q^{2}-1)}{2}z^{q-1}+(q^{2}-1)z^{q}+\frac{q(q-1)^{2}}{2}z^{q+1}.
		\end{align}
	\end{Lemma}
	
	We now use the geometry of ovals in $\PG(2,q)$ to show that the minimum-weight codewords of	the associated MDS codes form $q$-ary	$2$-designs when $q$ is even.
	
	\begin{Theorem}\label{oval-q-design}
		Let $q=2^m$ with $m\ge2$, let $\mathcal{O}$ be an oval in $\PG(2,q)$, and let $\C_{\mathcal{O}}$ be the associated $[q+1,3,q-1]_q$ code. Then $\Aw_{q-1}(\C_{\mathcal{O}})$ forms a $q$-ary $2$-$(q+1,q-1,(q-2)/2)$ design.
	\end{Theorem}
	\begin{proof}
		Let $\mathcal{O}=\{P_1,\ldots,P_{q+1}\}$ and write $G=G_{\mathcal{O}}$. For distinct points $P_s,P_t\in\mathcal{O}$, let $\ell_{st}=\overline{P_sP_t}$ denote their secant. By Lemma~\ref{Ha}, a nonzero codeword $\bm{c}=\bm{x}G$ has weight $q-1$ if and only if its associated line $\ell_{\bm{x}}$ is a secant of $\mathcal{O}$.
		
		Fix $1\le i<j\le q+1$ and $a,b\in\F_q^*$. Consider the projective point $R=\left[\bm{p}_j-\frac{b}{a}\bm{p}_i\right]$.
		It lies on $\ell_{ij}$, is distinct from $P_i,P_j$, and does not belong to the hyperoval $\mathcal{H}=\mathcal{O}\cup\{K\}$, where $K$ is the nucleus of $\mathcal{O}$. Indeed, $\ell_{ij}$ meets $\mathcal{H}$ exactly in $P_i$ and $P_j$.
		Every line through $R$ meets $\mathcal{H}$ in zero or two points. The $q+2$ points of $\mathcal{H}$ are therefore partitioned into pairs by the lines through $R$, so exactly $(q+2)/2$ of those lines meet $\mathcal{H}$ in two points. Exactly one of these lines, $\overline{RK}$, passes through the nucleus and is tangent to $\mathcal{O}$. The other $q/2$ lines are secants of $\mathcal{O}$. One is $\ell_{ij}$, leaving $(q-2)/2$ secants through $R$ distinct from $\ell_{ij}$.
		
		We show that these $(q-2)/2$ secants are in bijection with the codewords $\bm{c}\in\Aw_{q-1}(\C_{\mathcal{O}})$ satisfying $c_i=a$ and $c_j=b$. Let $\ell_{st}\ne\ell_{ij}$ be such a secant. Since $P_{i}\notin\ell_{st}$, the determinant $\det(\bm{p}_s,\bm{p}_t,\bm{p}_i)\ne 0$. Define the nonzero row vector
		\[
		\bm{x}=\alpha(\bm{p}_s\times\bm{p}_t),\qquad \alpha=\frac{a}{\det(\bm{p}_s,\bm{p}_t,\bm{p}_i)}.
		\]
		Since $R\in\ell_{st}$, $\det(\bm{p}_s,\bm{p}_t, \bm{p}_j-\frac{b}{a}\bm{p}_i)=0$. Then
		\[
		\det(\bm{p}_s,\bm{p}_t,\bm{p}_j)=\frac{b}{a}\det(\bm{p}_s,\bm{p}_t,\bm{p}_i).
		\]
		It follows that $\bm{c}=\bm{x}G$ has $c_i=a$ and $c_j=b$, and its only zero coordinates are $s,t$. Conversely, if a minimum-weight codeword satisfies these two coordinate conditions and has zeros at $s,t$, then the same determinant identity gives $R\in\ell_{st}$. Since $a,b\ne0$, the line $\ell_{st}$ differs from $\ell_{ij}$. The construction is unique in both directions, so
		\[
		|\{\bm{c}\in\Aw_{q-1}(\C_{\mathcal{O}}):c_i=a, c_j=b\}|=\frac{q-2}{2}.
		\]
		As this number is independent of $i,j,a,b$, the required $q$-ary $2$-design property follows.
	\end{proof}
	
	Combining the preceding results, we obtain a complete characterization of when the codewords of every nonzero weight in an $[q+1,3,q-1]_q$ MDS code form $q$-ary $2$-designs.

	\begin{Corollary}\label{every-k3MDS}
		Let $q\ge 3$ be a prime power, and let $\C$ be any	$[q+1,3,q-1]_q$ MDS code. Then the following statements are equivalent:
		\begin{enumerate}[(1)]
			\item $q$ is even;
			\item $\Aw_w(\C)$ forms a $q$-ary $2$-design for some $w\in\{q-1,q,q+1\}$;
			\item $\Aw_w(\C)$ forms a $q$-ary $2$-design for every $w\in\{q-1,q,q+1\}$.
		\end{enumerate}
		Moreover, when these conditions hold, the indices of the $q$-ary $2$-designs are $(q-2)/2$, $1$, and $q/2$, respectively.
	\end{Corollary}
	\begin{proof}
		By Corollary~\ref{cor:MDS-some-all}, statements (2) and (3) are equivalent.
		Suppose that (3) holds. Then $\Aw_{q-1}(\C)$ forms a $q$-ary $2$-design.
		By Lemma~\ref{lem:oval-weight-distribution}	and \eqref{q-lambda-i}, its index is $(q-2)/2$.
		Since this index is an integer, $q$ must be even.
		Conversely, suppose that $q$ is even. By the MDS--oval correspondence, $\C$ is monomially equivalent to an oval code.
		Theorems~\ref{oval-q-design} and~\ref{monomial-qtD} then imply that	$\Aw_{q-1}(\C)$ forms a $q$-ary $2$-design.
		Thus, (2) holds. Therefore, the three statements are equivalent. The remaining indices follow similarly from \eqref{q-lambda-i} and the	weight distribution in Lemma~\ref{lem:oval-weight-distribution}.
	\end{proof}

	The characterization also applies to the dual codes, by the duality principle established in Section 4.
	\begin{Corollary}\label{cor:dual-oval-complete}
		Let $q\ge3$ be a prime power, and let $\C$ be any	$[q+1,q-2,4]_q$ MDS code. Then $\Aw_w(\C)$ forms a $q$-ary $2$-design for every nonzero weight $w\in W(\C)$ if and only if $q$ is even. In particular, when $q$ is even, $\Aw_4(\C)$ forms a $q$-ary $2$-$(q+1,4,(q-2)/2)$ design.
	\end{Corollary}
	
	\begin{proof}
		Since $\C^{\perp}$ is an $[q+1,3,q-1]_q$ MDS code, we have	$\min\{d(\C),d(\C^{\perp})\}\ge2$.
		By Theorem~\ref{CCperp-qt} with $t=2$, $\Aw_w(\C)$ forms a $q$-ary $2$-design for every nonzero weight $w$ if and only if the same	holds for $\C^{\perp}$.
		Corollary~\ref{every-k3MDS} then establishes the claimed equivalence.
		Moreover, the weight distribution of an	MDS code gives	$A_4=(q-1)\binom{q+1}{4}$.
		Applying \eqref{q-lambda-i} yields $\lambda=(q-2)/2$, as desired.
	\end{proof}

	Although the codewords of every nonzero weight form $q$-ary $2$-designs for an $[q+1,3,q-1]_{q}$ code with even $q$, their design strength cannot be increased.
	
	\begin{Proposition}\label{prop:oval-not-three}
		Let $q=2^m$ with $m\ge2$, and let $\C$ be any $[q+1,3,q-1]_q$ MDS code. For every $w\in W(\C)$,	the set $\Aw_w(\C)$ does not form a $q$-ary $3$-design.
	\end{Proposition}
	\begin{proof}
		If $\Aw_w(\C)$ formed a $q$-ary $3$-design with index $\lambda_w$, then its index is given by $\lambda_w=\frac{A_w\binom{w}{3}}	{(q-1)^3\binom{q+1}{3}}$.
		Substituting the weight distribution from \eqref{eq:oval-weight-distribution}, we obtain
		\begin{align*}
			\lambda_{q-1}=\frac{(q-2)(q-3)}{2(q-1)^2},\quad	\lambda_q=\frac{q-2}{(q-1)^2},\quad	\lambda_{q+1}=\frac{q}{2(q-1)}.
		\end{align*}
		Since $q\ge4$, each of these numbers is strictly less than one and therefore cannot be the index of a $q$-ary $3$-design.
	\end{proof}

	As a special case of \cite[Theorem~9]{WuITW}, the codewords of every nonzero weight in a doubly-extended Reed--Solomon (DRS) code with parameters $[q+1,3,q-1]_q$ form $q$-ary $2$-designs when $q$ is even. 
	Unlike the corresponding result for DRS codes,
	Corollary~\ref{every-k3MDS} applies to all $[q+1,3,q-1]_q$ MDS codes, without	requiring monomial equivalence to a DRS code.
	
	\begin{Example}
		Let $q=8$. The enumeration in \cite[Table~1]{Beelen2017GRS}	implies the existence of $[9,3,7]_8$ MDS codes that are not	monomially equivalent to any DRS code.
		Let $\C$ be such a code. By Corollary~\ref{every-k3MDS}, the sets $\Aw_7(\C)$, $\Aw_8(\C)$, and $\Aw_9(\C)$ form $8$-ary $2$-$(9,7,3)$, $2$-$(9,8,1)$,
		and $2$-$(9,9,4)$ designs, respectively.
	\end{Example}

	\section{$q$-Ary $2$-Designs from BCH Codes and Their Duals}
	\setcounter{equation}{0}
	
	BCH codes are an important class of linear codes with rich algebraic structure. 
	In this section, we investigate $q$-ary $2$-designs arising from narrow-sense antiprimitive	BCH codes of length $q+1$ and their dual trace codes.
	We first extend the $3$-transitivity result in	\cite{3-tran-BCH} to a broader range of parameters.
	We then show that the codewords of every nonzero weight in certain trace codes form $q$-ary $2$-designs.
	The duality result from Section~4 gives the same conclusion	for the corresponding BCH codes.	
	Finally, we determine cases where the minimum-weight codewords do not form a $q$-ary $2$-design.
	For the BCH codes $\C_{(q,q+1,4,1)}$ with $q=2^{m}$, these results yield a complete characterization in terms of the parity of $m$.

	Let $p$ be a prime, $q=p^{m}$ and $\delta=p^{r}$ where $1\le r<m$. 
	Define $U=\{u\in \mathbb{F}_{q^{2}}^{*} : u^{q+1}=1\}$, which is the cyclic subgroup of $\mathbb{F}^{*}_{q^{2}}$ of order $q+1$. 
	Fix a generator $\gamma$ of $U$.
	The BCH code $\C_{(q,q+1,\delta,1)}$ is given by
	\begin{equation}\label{BCH-code}
		\C_{(q,q+1,\delta,1)}
		=\left\{(c_u)_{u\in U}\in\mathbb F_q^{q+1}: \sum_{u\in U}c_u u^i=0,\ 1\le i\le \delta-1\right\}.
	\end{equation}
	Let $\Tr_{q^{2}/q}(x)=x+x^{q}$ denote the trace function from $\mathbb{F}_{q^{2}}$ to $\mathbb{F}_{q}$. Define the trace code
	\begin{equation}
		\C_{\{1,\dots,\delta-1\}}
		=\left\{\left(\Tr_{q^2/q}\left(\sum_{i=1}^{\delta-1}a_i u^i\right)\right)_{u\in U}:	a_i\in\mathbb F_{q^2}\right\}.
	\end{equation}
	By \cite{SunDing2024}, these codes are dual to each other and have   parameters $[q+1,q+3-2\delta]_{q}$ and $[q+1,2\delta-2]_{q}$, respectively.
	Moreover, their minimum distances satisfy 
	\begin{align*}
		\delta\le d(\mathcal C_{(q,q+1,\delta,1)})\le 2\delta-1,\\
		q-2\delta+3\le d(\mathcal C_{\{1,\ldots,\delta-1\}})\le q-2\delta+4.
	\end{align*}
	
	We next recall the relevant projective linear group action.
	An element of $\PGL(2,q)$ represented by $\begin{pmatrix}a&b\\c&d\end{pmatrix}\in\mathrm{GL}(2,q)$
	acts on $\PG(1,q)$ by
	\[
	(x_{0}:x_{1})\mapsto (ax_{0}+bx_{1} : cx_{0}+dx_{1}).
	\]
	Under the usual identification $\PG(1,q)=\mathbb{F}_{q}\cup\{\infty\}$, this action can also be expressed as the linear fractional transformation
	$x\mapsto\frac{ax+b}{cx+d}$
	where $-d/c$ is mapped to $\infty$, and $\infty$ is mapped to $a/c$ if $c\ne 0$, or to $\infty$ if $c=0$.

	Let $\Stab_{U}=\{g\in\PGL(2,q^{2}) : g(U)=U\}$ be the setwise stabilizer of $U$ under the action of $\PGL(2,q^{2})$ on $\PG(1,q^{2})$.
	The structure of this stabilizer is described in the following lemma.	
	\begin{Lemma}\label{Stab-U-3t}
		\cite[Proposition 6]{3-tran-BCH} The setwise stabilizer $\Stab_{U}$ can be expressed as
		\[
		\left \{
		\begin{pmatrix}
			\beta_{2}^{q}&\beta_{1}^{q} \\
			\beta_{1}&\beta_{2}
		\end{pmatrix}\in \PGL(2,q^{2}) : \beta_{1}^{q+1}\ne \beta_{2}^{q+1}	
		\right \}.
		\]
		Moreover, $\mathrm{Stab}_{U}$ is conjugate in $\PGL(2,q^2)$ to $\PGL(2,q)$, and its action on $U$ is equivalent to the triply transitive action of $\PGL(2,q)$ on $\mathrm{PG}(1,q)$.
	\end{Lemma}
	
	To exploit this triple transitivity, we first identify a suitable lifted group acting on the trace code.
	Define
	\begin{align}\label{eq:F}
		\mathcal{F}:=\left\{
		\Tr_{q^2/q}\left(\sum_{i=1}^{\delta-1}a_i u^i\right)  \in\mathbb{F}_{q^{2}}[u]/\langle u^{q+1}-1\rangle :	a_i\in\mathbb F_{q^2}
		\right\}.
	\end{align}
	Thus $\C_{\{1,\dots,\delta-1\}}=\{(f(u))_{u\in U}:f\in \mathcal{F}\}$. Define the lift of Stab$_{U}$ by
	\[
	\overline{\Stab}_{U} = \left \{ \begin{pmatrix}
		\beta_{2}^{q}&\beta_{1}^{q} \\
		\beta_{1}&\beta_{2}
	\end{pmatrix}\in \text{GL}(2,q^{2}) : \beta_{1},\beta_{2}\in\mathbb{F}_{q^{2}},\beta_{1}^{q+1}\ne \beta_{2}^{q+1} \right \}.
	\]
	For $M=\begin{pmatrix}
		\alpha_{1}&\alpha_{2} \\
		\beta_{1}&\beta_{2}
	\end{pmatrix}\in \overline{\text{Stab}}_{U}$ and $\overline{G}=M^{-1}$, define
	\[
	(\overline{G}\circ f)(u) = (\beta_{1}u+\beta_{2})^{(\delta-1)(q+1)}f\left (\frac{\alpha_{1}u+\alpha_{2}}{\beta_{1}u+\beta_{2}}\right ), \quad f\in \mathcal{F}.
	\]
	This expression is well defined on $U$. Indeed,  $\beta_1u+\beta_2=0$ for some $u\in U$ implies that $\beta_1^{q+1}=\beta_2^{q+1}$.

	The following lemma shows that these transformations preserve $\mathcal F$ and define a group action.
	
	\begin{Lemma}\cite[Lemma 10]{3-tran-BCH} \label{TF-representation}
		Let $p$ be prime, $q=p^{m}$ and $\delta=p^{r}$ with $1\le r<m$. Let $E$ denote the $2\times 2$ identity matrix. 
		For any $\overline{G}_{1},\overline{G}_{2}\in \overline{\text{Stab}}_{U}$ and any $f_{1},f_{2}\in\mathcal{F}$, the following hold:
		\begin{enumerate}
			\item[(1)] $\overline{G}_{1}\circ f_{1}\in \mathcal{F}$;
			\item[(2)] $E\circ f_{1} = f_{1}$;
			\item[(3)] $(\overline{G}_1\overline{G}_2)\circ f_{1} = \overline{G}_{1}\circ(\overline{G}_{2}\circ f_{1})$;
			\item[(4)] $\overline{G}_{1}\circ (af_{1}+bf_{2}) = a(\overline{G}_{1}\circ f_{1}) + b(\overline{G}_{1}\circ f_{2})$ for all $a,b\in\mathbb{F}_{q}$.
		\end{enumerate}
	\end{Lemma}
	
	In \cite[Theorem~21]{3-tran-BCH}, the $3$-transitivity	result was established under the assumption	$q=\delta^s$. However, the preceding lemma applies to the more general setting $q=p^m$ and $\delta=p^r$ with $1\le r<m$.
	\begin{Theorem}\label{C13-3transitive}
		Let $p$ be prime, $q=p^{m}$ and $\delta=p^{r}$ with $1\le r<m$.
		Define
		\[
		G = \left\{ \left( ((\beta_{1}u+\beta_{2})^{(\delta-1)(q+1)})_{u\in U}; 
		\begin{pmatrix}
			\beta_{2}^{q}&\beta_{1}^{q} \\
			\beta_{1}&\beta_{2}
		\end{pmatrix}^{-1}\right) : \beta_{1},\beta_{2}\in\mathbb{F}_{q^{2}},\beta_{1}^{q+1}\ne\beta_{2}^{q+1}\right\}.
		\]
		Then $G\le\MAut(\C_{\{1,\dots,\delta-1\}})$. 
		In particular, $\Aut(\mathcal C_{\{1,\dots,\delta-1\}})$ is $3$-transitive on  $U$.
	\end{Theorem}
	\begin{proof}
		For each $M\in\overline{\Stab}_U$, the quantity $\beta_1u+\beta_2$ is nonzero for every $u\in U$. Hence  $(\beta_{1}u+\beta_{2})^{(\delta-1)(q+1)}\in\F_{q}^{*}$.
		Then the coordinate multipliers shown by the first components of elements in $G$ are nonzero elements of $\mathbb F_q$.
		By Lemma~\ref{TF-representation}, the
		corresponding transformations preserve the trace polynomial space in~\eqref{eq:F} and define a group action.
		Consequently, $G\le\MAut(\mathcal C_{\{1,\ldots,\delta-1\}})$.
		By Lemma~\ref{Stab-U-3t}, the permutation part of $G$ acts $3$-transitively on $U$.
		The result follows.
	\end{proof}
	This extension also gives ordinary $3$-designs in the corresponding BCH codes and their duals.
	\begin{Proposition}
		Let $p$ be prime, $q=p^{m}$ and $\delta=p^{r}$ with $1\le r<m$.
		For every $w\in W(\mathcal C_{\{1,\ldots,\delta-1\}})$ with $w\ge3$, the set $\mathcal B_w(\mathcal C_{\{1,\ldots,\delta-1\}})$ forms a $3$-design. The same holds for $\mathcal B_w(\mathcal C_{(q,q+1,\delta,1)})$ whenever $w\in W(\mathcal C_{(q,q+1,\delta,1)})$ and $w\ge3$.
	\end{Proposition}
	\begin{proof}
		By Theorem~\ref{C13-3transitive}, the automorphism group of the trace code $\mathcal C_{\{1,\ldots,\delta-1\}}$ is $3$-transitive. Hence its codewords of every nonzero weight $w\ge3$ support a $3$-design by Theorem~\ref{transitive-t-design}. Since $\C_{(q,q+1,\delta,1)}$ and $\C_{\{1,\dots,\delta-1\}}$  are dual, their monomial automorphism groups have the same induced permutation group on $U$. The BCH code therefore also has a $3$-transitive automorphism group, and the same conclusion follows.
	\end{proof}
	
	From now on, we assume $q=2^{m}$ and $\delta=2^{r}$, where $1\le r<m$. 
	For any $h\in\mathbb{F}_{q}^{*}$, we define the bijection $\varphi_{h} : U\to U$ by 
	\begin{align*}
		\varphi_{h}(x) =  \frac{B^{q}x+A^{q}}{Ax+B}, \qquad A = h+1, B=\gamma^{-1}+\gamma h. 
	\end{align*}
	In characteristic $2$,
	\[
	B^{q+1}-A^{q+1}=h(\gamma^2+\gamma^{-2})=h(\gamma+\gamma^{-1})^2\ne0,
	\]
	since $h\ne0$ and $\gamma\ne\gamma^{-1}$. Theorem~\ref{C13-3transitive} implies that the monomial transformation
	\begin{align}\label{eq:gh}
		g_{h}&=\left( \left (a_{h}(u)\right )_{u\in U}; \varphi_{h}^{-1}\right)\in\MAut(\C_{\{1,\dots,\delta-1\}}),
	\end{align}
	where $a_{h}(u)=(Au+B)^{(\delta-1)(q+1)}$. Indeed $\varphi_{h}$ is represented by $\begin{pmatrix}B^{q}&A^{q} \\	A&B	\end{pmatrix}$, and hence
	\begin{align*}
		g_{h}=\left( \left ((Au+B)^{(\delta-1)(q+1)}\right )_{u\in U}; 
		\begin{pmatrix}	B^{q}&A^{q} \\	A&B	\end{pmatrix}^{-1}\right).
	\end{align*}
	Consequently, for every $\bm{c} = (c_{u})_{u\in U}\in \C_{\{1,\dots,\delta-1\}}$, 
	\begin{align}\label{eq:ghc}
		g_{h}(\bm{c}) = \left (a_{h}(u)c_{\varphi_{h}(u)}\right )_{u\in U}\in\C_{\{1,\dots,\delta-1\}}.
	\end{align}
	By the definition of $\varphi_{h}$, we have $\varphi_{h}(\gamma)=\gamma$ and  $\varphi_{h}(\gamma^{-1})=\gamma^{-1}$. 
	Thus $\varphi_{h}$ fixes $\gamma$ and $\gamma^{-1}$. 
	Moreover, the multipliers at these two fixed coordinates satisfy 
	\begin{align}\label{eq:h-6}
		\frac{a_{h}(\gamma)}{a_{h}(\gamma^{-1})} = h^{-(\delta-1)(q+1)} = h^{-2(\delta-1)}.
	\end{align}
	We are now ready to establish the desired $q$-ary $2$-design property in $\C_{\{1,\dots,\delta-1\}}$.

	\begin{Theorem}\label{C123Q-2}
		Let $q=2^{m}$ and $\delta=2^{r}$, where $1\le r<m$ and $\gcd(r,m)=1$. Let $\C=\C_{\{1,\dots,\delta-1\}}$. Then $\mathcal{A}_{w}(\C)$ forms a $q$-ary $2$-$(q+1,w,\lambda)$ design for every nonzero weight $w\in W(\C)$, where $\lambda = \frac{A_{w}(\C)w(w-1)}{(q-1)^{2}q(q+1)} $.
	\end{Theorem}
	\begin{proof}
		For a vector $\bm{y}\in\mathbb{F}_{q}^{U}$ of weight $2$, define
		\[
		\lambda_{w,\bm{y}}=\left|\left\{\bm{c}\in\mathcal{A}_{w}(\C):	\bm{c}\text{ covers }\bm{y}\right\}\right|.
		\]
		It suffices to prove that $\lambda_{w,\bm{y}}$ is independent of the choice of $\bm{y}$.
		
		Let $\bm{y}$ be an arbitrary vector of weight $2$, with
		$\supp(\bm{y})=\{u_1,u_2\}.$
		By Theorem~\ref{C13-3transitive} and Theorem~\ref{transitive-design},
		it is sufficient to consider vectors with support $\{\gamma,\gamma^{-1}\}$.
		For $z\in\mathbb{F}_{q}^{*}$, define the vector $\bm{y}_{z}$ by $(\bm{y}_{z})_{\gamma}=z,\ (\bm{y}_{z})_{\gamma^{-1}}=1$,
		and zero elsewhere.
		Since $\C$ is linear, every vector of weight $2$ with support $\{\gamma,\gamma^{-1}\}$ is equivalent, with respect to the covering number $\lambda_{w,\bm{y}}$, to some $\bm{y}_{z}$.
		
		Let $z_{1},z_{2}\in\mathbb{F}_{q}^{*}$. Observe that 
		\[
		\gcd(2(\delta-1),q-1)=\gcd(2^{r}-1,2^{m}-1)=2^{\gcd(r,m)}-1=1.
		\]
		Therefore the map $h\longmapsto h^{2(\delta-1)}$ is a bijection. We may thus choose $h\in\mathbb{F}_{q}^{*}$ such that $ h^{2(\delta-1)}=z_{1}z_{2}^{-1}$.
		Consider the automorphism $g_h$ in \eqref{eq:gh}. Since $\varphi_h$ fixes both $\gamma$ and $\gamma^{-1}$, the support of $g_h(\bm{y}_{z_1})$ is again $\{\gamma,\gamma^{-1}\}$.
		By \eqref{eq:ghc} and \eqref{eq:h-6},
		\begin{align*}
			\frac{(g_h(\bm{y}_{z_1}))_{\gamma}}{(g_h(\bm{y}_{z_1}))_{\gamma^{-1}}}=\frac{a_h(\gamma)z_1}{a_h(\gamma^{-1})}
			=h^{-2(\delta-1)}z_1=z_2.
		\end{align*}
		Hence $\eta g_h(\bm{y}_{z_1})=\bm{y}_{z_2}$ for some $\eta\in\mathbb{F}_{q}^{*}$.
		Since both $g_h$ and multiplication by $\eta$ are monomial automorphisms of $\C$, they induce a bijection between the codewords of weight $w$ covering $\bm{y}_{z_1}$ and those covering $\bm{y}_{z_2}$. Therefore $\lambda_{w,\bm{y}_{z_1}}= \lambda_{w,\bm{y}_{z_2}}$.
		Combining the above arguments shows that $\lambda_{w,\bm{y}}$
		is constant over all vectors $\bm{y}$ of weight $2$. Hence $\mathcal{A}_{w}(\C)$ forms a $q$-ary $2$-design.
		The index $\lambda=\lambda_{2}$ follows from \eqref{q-lambda-i}.
	\end{proof}

	Combining the preceding result with	the duality principle yields the following theorem.
	\begin{Theorem}\label{trace-BCH}
		Let $q=2^{m}$ and $\delta=2^{r}$, where $1\le r<m$ and $\gcd(r,m)=1$. 
		Then, for every $w\in W(\mathcal C_{\{1,\ldots,\delta-1\}})$, $\mathcal A_w(\mathcal C_{\{1,\ldots,\delta-1\}})$ forms a $q$-ary $2$-design. Likewise, for every $w\in W(\mathcal C_{(q,q+1,\delta,1)})$, $\mathcal A_w(\mathcal C_{(q,q+1,\delta,1)})$ forms a $q$-ary $2$-design.
	\end{Theorem}
	
	\begin{proof}
		The assertion for $\mathcal C_{\{1,\dots,\delta-1\}}$ follows from Theorem~\ref{C123Q-2}. By \cite{SunDing2024}, $\mathcal C_{(q,q+1,\delta,1)}=\mathcal C_{\{1,\dots,\delta-1\}}^\perp$.
		Moreover, since $r<m$, we have $2\delta\le q$, and hence
		\[
		d(\mathcal C_{\{1,\dots,\delta-1\}})\ge q-2\delta+3\ge3,
		\qquad
		d(\mathcal C_{(q,q+1,\delta,1)})\ge\delta+1\ge3.
		\]
		Theorem~\ref{CCperp-qt} therefore applies with $t=2$ and gives the result for the BCH code.
	\end{proof}
	
	We now turn to the special case $\delta=4$ and determine some indices of corresponding $q$-ary $2$-designs. 
	
	\begin{Proposition}\label{prop:BCH-trace-indices}
		Let $q=2^m$, where $m\ge5$ is odd, and let $\C=\C_{(q,q+1,4,1)}$ and $\mathcal{D}=\C_{\{1,2,3\}}$.
		Then the following hold:
		\begin{enumerate}
			\item[(1)] For $w=6,7$, $\Aw_w(\C)$ forms a	$q$-ary $2$-$(q+1,w,\lambda_w)$ design,
			where $\lambda_6=(q-2)(q-8)/4!$, $\lambda_7=(q-2)(q-5)(q^2-14q+68)/5!$.
			\item[(2)] For $w=q-5,\ q-4$, $\Aw_w(\mathcal{D})$ forms a	$q$-ary $2$-$(q+1,w,\lambda_w)$ design,
			where $\lambda_{q-5}=(q-2)(q-5)(q-6)(q-8)/6!$, $\lambda_{q-4}=(q-2)(q-4)(q-5)/4!$.
		\end{enumerate}
	\end{Proposition}
	
	\begin{proof}
		Since $m$ is odd, it follows from Theorem~\ref{trace-BCH} that the codewords of every nonzero weight in	both $\C$ and $\mathcal{D}$ form $q$-ary $2$-designs. It remains to determine their indices at the specified weights.
		
		By \cite[Theorems~34 and 37]{BCH-4-design},	$\C$ and $\mathcal{D}$ have parameters
		$[q+1,q-5,6]_q$ and $[q+1,6,q-5]_q$, respectively.
		For a $q$-ary $2$-$(n,w,\lambda_{w})$ design with block set $\Aw_{w}$, we have $\lambda_{w}=\frac{|\Aw_{w}|\binom{w}{2}}{(q-1)^2\binom{n}{2}}$.
		By \cite{BCH-4-design,Yan}, the numbers of codewords of weights $6$ and $7$ in $\C$
		satisfy $A_6(\C)=\frac{(q-8)(q-1)}{30}\binom{q+1}{4}$ and $A_7(\C)=(q-1)\binom{q+1}{7}-(q-5)A_6(\C)$.		
		Substituting these expressions into $\lambda_{w}$ gives $\lambda_6$ and $\lambda_7$. This proves (1).
		The indices in (2) follow from	\cite[Theorems~5.12 and 5.13]{wu}.
	\end{proof}
	
	To end this section, we give a complete characterization of the minimum-weight codewords of $\C_{(q,q+1,4,1)}$ in terms of $q$-ary designs.

	\begin{Theorem}\label{th:BCH4-iff-modd}
		Let $q=2^{m}$, where $m\ge 4$.
		Let $\C=\C_{(q,q+1,4,1)}$ be the code with minimum distance $d$. Then $\Aw_{d}(\C)$ forms a $q$-ary $2$-design if and only if $m$ is odd.
	\end{Theorem}
	\begin{proof}
		By~\cite[Theorem 34]{BCH-4-design}, $\C$ has parameters $[q+1,q-5,d]_{q}$, where $d=6$ if $m$ is odd and $d=5$ if $m$ is even.
		Further if $m$ is even, then the BCH code $\C$ supports a $3$-$(q+1,5,1)$ design at the minimum weight. Hence $A_{5}=q(q+1)(q-1)^{2}/60 $.
		If $\Aw_{5}(\C) $ is a $q$-ary $2$-design, then $A_{5}= \lambda(q+1)q(q-1)^{2}/20$, giving a contradiction that $\lambda$ is not an integer.
		If $m$ is odd, it follows from Theorem~\ref{trace-BCH} that $\Aw_{d}(\C)$ forms a $q$-ary $2$-design. This completes the proof.
	\end{proof}

	\section{Summary and Concluding Remarks}
	\setcounter{equation}{0}
	
	This paper established several criteria for determining	when the codewords of every nonzero weight in a linear code form a $q$-ary $t$-design. These results extend the study of $q$-ary designs from certain smaller weights	to all nonzero weights under suitable conditions. We first showed that this property is preserved under monomial equivalence and examined its relationship with $t$-regularity of the code. In particular, when $t\le d< 2t$, $t$-regularity together with the $q$-ary $t$-design property at the minimum weight guarantees the same property at every nonzero weight (Theorem~\ref{t-regular+d->all}). A second criterion reduces the verification to weights $d\le w\le s^{\perp}+t-1$ under the hypotheses of Theorem~\ref{th:s'+t-1}. For MDS codes, we further proved that the $q$-ary $2$-design property at every nonzero weight is equivalent to that at any fixed nonzero weight (Corollary~\ref{cor:MDS-some-all}).
	
	Considering also the $q$-ary designs in the dual code, we established a duality principle using additive characters and Krawtchouk polynomials: if $1\le t\le\min\{d,d^{\perp}\}$, the codewords of every nonzero weight in a linear code form a $q$-ary $t$-design if and only if it holds for its dual (Theorem~\ref{CCperp-qt}).

	As an application, we characterized the $[q+1,3,q-1]_q$	MDS codes in which the codewords of every nonzero weight form $q$-ary $2$-designs. This occurs precisely when	$q$ is even (Corollary~\ref{every-k3MDS}).
	The same characterization holds for their $[q+1,q-2,4]_q$ dual codes	(Corollary~\ref{cor:dual-oval-complete}).
	However, for even $q$, the codewords of any nonzero	weight in an $[q+1,3,q-1]_q$ MDS code do not form a $q$-ary $3$-design	(Proposition~\ref{prop:oval-not-three}).
	We also constructed infinite families of $q$-ary $2$-designs from narrow-sense antiprimitive BCH codes and their duals.
	For $q=2^m$ and $\delta=2^r$ with $1\le r<m$ and $\gcd(r,m)=1$,	the codewords of every nonzero weight in $\C_{(q,q+1,\delta,1)}$ and its dual form $q$-ary $2$-designs	(Theorem~\ref{trace-BCH}).
	Furthermore, when $\delta=4$ and $m\ge4$, we obtained a complete
	characterization: the minimum-weight codewords of $\mathcal C_{(q,q+1,4,1)}$ form a $q$-ary $2$-design if and only if $m$ is odd
	(Theorem~\ref{th:BCH4-iff-modd}).

	These results suggest several directions for further investigation:
	\begin{enumerate}[(1)]
		\item Can the verification range $d\le w\le s^{\perp}+t-1$ in Theorem~\ref{th:s'+t-1} be shortened? More generally, which subsets of nonzero weights suffice to ensure that the codewords of every nonzero weight form $q$-ary $t$-designs?		
		For what kinds of codes can this verification be reduced to the minimum weight, except for the MDS codes in the $t=2$ case?
		
		\item Let $q=2^m$ and $\delta=2^r$, where $1\le r<m$ and $\gcd(r,m)>1$.	For which nonzero weights do the codewords in	$\mathcal C_{(q,q+1,\delta,1)}$ and its dual form $q$-ary $2$-designs?
		
		\item Let $q=2^m$ and $\delta=2^r$, where $1\le r<m$ and $\gcd(r,m)=1$.	Determine whether $\C_{(q,q+1,\delta,1)}$ supports ordinary $4$-designs for $\delta\ne4$.
		
	\end{enumerate}
	
	\section*{Acknowledgment}
	This work was supported by National Natural Science Foundation of China under Grant 12571346. 	
	The authors used OpenAI ChatGPT as an auxiliary tool for language polishing and presentation improvement and to assist with supplementary consistency checks. The research ideas, mathematical derivations, proofs, and conclusions were developed by the authors. All suggestions from the tool were critically reviewed and revised as appropriate. The authors independently verified the mathematical content and references and assume full responsibility for the accuracy and integrity of the final manuscript.

\end{document}